\documentclass[review]{siamart220329}

\usepackage{amsfonts}
\usepackage{graphicx}
\usepackage{epstopdf}
\usepackage{algorithmic}

\ifpdf
  \DeclareGraphicsExtensions{.eps,.pdf,.png,.jpg}
\else
  \DeclareGraphicsExtensions{.eps}
\fi

\usepackage{enumitem}
\setlist[enumerate]{leftmargin=.5in}
\setlist[itemize]{leftmargin=.5in}
\newsiamremark{remark}{Remark}
\newsiamremark{hypothesis}{Hypothesis}
\crefname{hypothesis}{Hypothesis}{Hypotheses}
\newsiamthm{claim}{Claim}

\usepackage{amsopn}

\newcommand{\TheTitle}{%
  Localized subspace iteration methods for elliptic multiscale problems
}
\newcommand{\TheShortTitle}{%
Localized Subspace Iteration 
}
\newcommand{\TheShortName}{%
  Xiaofei Guan,
  Lijian Jiang,
  Yajun Wang,
  Zihao Yang
}

\newcommand{\TheTJmathAddress}{%
  School of Mathematical Sciences, and Key Laboratory of Intelligent Computing and Applications (Ministry of Education),  Tongji University, Shanghai 200092, China (\email{guanxf@tongji.edu.cn,~ljjiang@tongji.edu.cn,~1910733@tongji.edu.cn}).
}
\newcommand{\TheNPUmathAddress}{%
  School of Mathematics and Statistics, Northwestern Polytechnical University, Xi'an 710072, China (\email{yangzihao@nwpu.edu.cn}).
}

\newcommand{\TheFunding}{%
  Submitted to the editors DATE.
  \funding{The work of the first author was supported by the National Science Foundation of China under grant 12271409, the Interdisciplinary Project in Ocean Research of Tongji University and the Fundamental Research Funds for the Central Universities. The work of the second author was supported by the National Science Foundation of China under grant 12271408.}
}

\author{
  Xiaofei Guan\thanks{\TheTJmathAddress}
  \and
  Lijian Jiang\footnotemark[2]
  \and
  Yajun Wang\footnotemark[2]
  \and
  Zihao Yang\thanks{\TheNPUmathAddress}
}

\graphicspath{{figures/}}
\usepackage{stmaryrd}
\usepackage{mathrsfs}
\usepackage{subfigure}
\usepackage{booktabs}
\usepackage{multirow}
\usepackage{amsopn}
\usepackage{cleveref}

\newtheorem{mythm}{Theorem}[section]

\newtheorem{mylem}{Lemma}[section]
\newtheorem{myassum}{Assumption}[section]

\newtheorem*{mypf}{Proof}
\newtheorem{myrema}{Remark}[section]

\title{{\TheTitle}\thanks{\TheFunding}}
\headers{\TheShortTitle}{\TheShortName}
\ifpdf%
  \hypersetup{%
    pdftitle={\TheTitle},
    pdfauthor={\TheShortName}
  }
\fi

\begin{document}

\maketitle
\begin{abstract}
This paper proposes localized subspace iteration (LSI) methods to construct generalized finite element basis functions for elliptic problems with multiscale coefficients.
The key components of the proposed method consist of the localization of the original differential operator and the subspace iteration of the corresponding local spectral problems, where the localization is conducted by enforcing the local homogeneous Dirichlet condition and the partition of the unity functions.
From a novel perspective, some multiscale methods can be regarded as one iteration step under approximating the eigenspace of the corresponding local spectral problems. Vice versa, new multiscale methods can be designed through subspaces of spectral problem algorithms.
Then, we propose the efficient localized standard subspace iteration (LSSI) method and the localized Krylov subspace iteration (LKSI) method based on the standard subspace and Krylov subspace, respectively.
Convergence analysis is carried out for the proposed method. Various numerical examples demonstrate the effectiveness of our methods. In addition, the proposed methods show significant superiority in treating long-channel cases over other well-known multiscale methods.
\end{abstract}


\section{Introduction}
 Multiscale problems are important problems in both scientific computing and engineering applications. Examples include diffusion in fractured media \cite{UR_DMfFM_1990} and deformation of composite materials with multiple nonseparated length scales \cite{M_MoDaAPiPM_1962}, etc.
However, the computation simulation by traditional numerical methods poses significant challenges because of the expensive computational cost for these problems \cite{IJ_CAFEMPAB_2000}.
This issue is caused by extremely fine computational grids that are required to resolve all relevant scales.
To this end, numerical homogenization \cite{AD_NHbLOD_2020,Houman_OAWFSaNH_2019} has been developed and replaced polynomial finite element ansatz functions with more general ansatz functions to overcome the global fine-scale computation.
The essence of these methods lies in designing particular finite element ansatz functions that can efficiently capture the problem's  multiscale information.

Babuska, et al.'s pioneering work on the Generalized Finite Element Method (GFEM) \cite{IJ_GFEM_1983} suggested that a specific form of basis function is desirable for one-dimensional problems with rough coefficients. Then the idea is extended to two-dimensional problems \cite{IGJ_SFEM_1994}.
The multiscale finite element method (MsFEM) \cite{HW_MsFEM_1997,Yalchin_MsFEM:TaA_2009} by Hou et al. is another significant milestone. They proposed to generate the harmonic extension of the standard Lagrangian finite element basis functions and to use these harmonic extensions afterward as ansatz functions in the Galerkin method.
This method has broad applicability, but it can lead to so-called resonance errors \cite{Hou_RtCRE_2004,Henning_OfMsFEM_2013} due to  the artificial boundary conditions of the local harmonic extension problems.
After that, an important method to avoid resonance errors was suggested in \cite{Houman_MBU_2007,AEJ_MMsLGI_2008}. It involves using known global information to figure out better boundary conditions for local problems.
If such global information is not available, more sophisticated constructions are necessary to derive suitable multiscale ansatz functions from local problems.

More recently, some new multiscale methods have emerged to tackle this challenge.
One way is to use local spectral problems to figure out what information is redundant and what information is important. Then, you can use this information to build multiscale ansatz functions, such as SGFEM \cite{Babuska_OLASfGFEM_2011,BLSS_MSGFFEM_2020,Ma_NDaAfGFEM_2022} and GMsFEM \cite{EGH_GMsFEM_2013}, etc.
This significantly improves computation accuracy.
Another method based on the orthogonal decomposition of the solution space was developed in \cite{AD_LOD_2014,HM_LODfBVP_2014}. It is capable of converting arbitrary finite element (FE) basis functions into ansatz functions incorporating multiscale information. The exponential decay of these multiscale ansatz functions allows them to be localized, which is also known as localized orthogonal decomposition (LOD).
By combining this technique with spectral problems, CEM-GMsFEM proposed in \cite{EYW_CEMGMs_2018,LCJ_CEMfP_2019}  is able to improve POD and GMsFEM.
There are also numerous methods for constructing some other multiscale  basis functions and their variations \cite{Allaire_Homogenization_1992,Hughes_VMS_1998,Wu_upscaling_2002,WeinanE_HMM_2003,Arbogast_mortar_2007,GJW_MMRfSEPEVS_2023}. The reference  list is  incomplete. Subsequently, one crucial question is: What are the fundamental principles governing the construction of these basis functions?
In our opinion, there are two fundamental principles for the construction of multiscale basis functions.
The first is that the basis functions must have localized support, and this ensures that the stiffness matrix is sparse. Furthermore, solving localized problems is computationally efficient and desirable for parallel techniques. This will significantly decrease the CPU time for building the basis functions.
Another fundamental principle is that local problems are connected to the inverse operator corresponding to the original problem. Several multiscale methods \cite{AD_LOD_2014,HM_LODfBVP_2014,EYW_CEMGMs_2018} use the inverse operator in advance, either locally or indirectly. By integrating these two principles, the local inverse operator is crucial in the design of multiscale basis functions.

An easy way to get multiscale basis functions from a standard set of finite element basis functions is to use the orthogonal decomposition technique in the LOD (\cite{AD_LOD_2014,HM_LODfBVP_2014}).
Employing the orthogonal decomposition technique several times will result in a more accurate basis function space sequence, and we refer to it as iterative orthogonal decomposition.
As proved in \cite{RPD2021NHBSS}, the trial and test function spaces obtained by the orthogonal decomposition technique are equivalent to the function spaces obtained by applying the inverse operator and the inverse conjugate operator to the initial basis function space, respectively.
The iterative function space sequence that results from iterative orthogonal decomposition converges to the eigenfunction subspace of the inverse operator.
At the same time, this iterative function space sequence is consistent with the subspace sequence $\{ A^k X_0 \}$ for solving the matrix eigenvalue problem $AX = \lambda X$ \cite{{Saad_Iterative_2003}}. These motivate us to construct multiscale basis function spaces using a variety of subspaces derived from corresponding spectral problems.


The basic idea of the LSI starts with designing the local inverse operators in each local domain $\omega_i$, where $\left\{\omega_i \right\}_i$ is an overlapping open cover of domain $\Omega$. Moreover, we opt to enforce a local homogeneous Dirichlet boundary condition on the operator to accomplish localization.
Then, novel multiscale methods can be developed by the subspace, which is initially utilized to deal with the spectral problem corresponding to the local inverse operator. Specifically, a localized standard subspace iteration (LSSI) method is proposed by utilizing standard subspaces \cite{Saad_Iterative_2003} of local inverse operators. Furthermore, a localized Krylov subspace iteration (LKSI) method is also developed by utilizing the Krylov subspaces \cite{Saad_Iterative_2003,Liesen_Krylov_2013} of local inverse operators. The proposed multiscale methods are applied to multiscale diffusion and elasticity equations, and they provide better stability than many other multiscale methods in dealing with long-channel problems.

The paper is organized as follows: In \Cref{section:num2}, an inherent relationship between the LOD and the subspace iteration is established. Then, in \Cref{section:num3}, two classic subspaces for solving spectral problems are introduced. \Cref{section:num4} presents two multiscale methods based on the subspace iteration. In \Cref{section:num5}, the convergence analysis of proposed methods is carried out. In \Cref{section:num6}, a few numerical examples are provided to demonstrate the efficacy and accuracy of the proposed methods. Subsequently, conclusions are made in \Cref{section:num7}.


\section{From LOD to spectral problem} \label{section:num2}

\subsection{Introduction to LOD}
Localized Orthogonal Decomposition (LOD) \cite{AD_LOD_2014,HM_LODfBVP_2014}, a highly effective method for constructing multiscale basis functions, has greatly driven the development of many multiscale models. Consider the following elliptic problem:
\begin{equation} \label{eq_govering}
	\mathcal{L} u(x)=f(x) \quad \text { in } \Omega,
\end{equation}
where $\Omega \in \mathbb{R}^d$ is a polyhedral Lipschitz domain, $f \in L^2(\Omega)$ denotes a given source term, and $\mathcal{L} $ is a linear partial differential operator with some high-contrast or high-oscillation multiscale coefficients. For simplicity, suppose the above equation satisfies the homogeneous Dirichlet boundary conditions $u = 0 \text { in } \partial \Omega $. Let $V$ denote a Sobolev space that match the problem, then the variational form of (\ref{eq_govering}) is given as
\begin{equation} \label{eq_variational_form}
	a(u, v)=\left( f ,v \right)_{L^2(\Omega)} \text { for all } v \in V,
\end{equation}
where $a(u, v)$ is a bounded sesquilinear form.
Let $\mathcal{T}_H$ denote a regular finite element mesh of the domain $\Omega$ into closed simplices. Let $\left\{l_j\right\}$ denote a set of finite element basis functions, such as the nodal basis functions (hat functions) of the Lagrange finite element space, and define $U_H = \text{span} \left\{ l_j \right\}$.

The LOD starts with a number of \textit{macroscopic quantities of interest}, which extract the desired information from the exact solution \cite{RPD2021NHBSS}. These continuous linear functionals are denoted as $q_j \in \mathcal{V}^*, j \in J$, where $J$ is the finite index set with $N:=|J|$.
Without loss of generality, we assume that these functionals are linearly independent. A canonical choice of functional $q_j$ is
\begin{equation}
	q_j:=\left(l_j, \bullet\right)_{L^2(\Omega)}.
\end{equation}
In fact, there are numerous alternative selections for the quantities of interest.
\begin{myrema}
    For another set of quantities of interest $ \tilde{q}_j \in \mathcal{V}^*$, it is necessary to make the assumption $ \tilde{q}_j \in \left[L^2(\Omega)\right]^*$.  According to Riesz representation theorem, there exists a unique $\tilde{l}_j \in \mathcal{V}$ such that
    \begin{equation}
        \tilde{q}_j(v) = \left(\tilde{l}_j, v \right)_{L^2(\Omega)}, \ \forall \  v \in V.
    \end{equation}
    Hence, each quantity of interest satisfying $ \tilde{q}_j \in \left[L^2(\Omega)\right]^*$ is uniquely associated with a corresponding basis function, and vice versa. It is worth noting that the $L^2$ inner product used here can be replaced by other well-defined inner products, such as weighted $L^2$ inner product defined by
    \begin{equation}
        \left(u, v \right)_{L^2(\Omega,\kappa)} := \int_{\Omega} \kappa u v.
    \end{equation}
\end{myrema}

Given the macroscopic quantities of interest $q_j$, we proceed to establish the kernel space
\begin{equation}
  W:=\left\{v \in V \mid q_j(v)=0 \text { for all } j \in J\right\}=\bigcap_{j \in J} \operatorname{ker} q_j.
\end{equation}
This space is sometimes referred to as the fine-scale space \cite{TGLJ_VMS_1998}, which encompasses fine-scale information that cannot be captured by the original basis function space $V_H$. It inspired adaptive methods involving the addition of basis functions \cite{MA_AVMS_2007}.
The LOD method's core idea is to find a function space that is orthogonal to the given kernel space $W$ in the sense of sesquilinear form $a\left(\bullet,\bullet \right)$, serving as the multiscale basis function space.
To this end, we define two projections $\mathcal{C}: {V} \rightarrow {W}$ and $\mathcal{C}^*: {V} \rightarrow {W}$, such that
\begin{equation}
  a(\mathcal{C} v, w)=a(v, w) \text{ and } a(w,\mathcal{C}^* v)=a(w, v) \ \forall \ v \in V, w \in W.
\end{equation}
A natural conclusion is that $\mathcal{C}= \mathcal{C}^*$ if $a\left(\bullet,\bullet \right)$ is Hermitian.
Using the kernel space $W$ and the projections $C$ and $C^*$, the trial and test spaces are constructed by
\begin{equation}
  \widetilde{U}_H:=(1-\mathcal{C}) V \quad \text { and } \quad \widetilde{V}_H:=\left(1-\mathcal{C}^*\right) V.
\end{equation}

\begin{mylem}
  Given the inf-sup condition, the spaces $\widetilde{U}_H$ and $\widetilde{V}_H$
  possess a dimension of $N:=|J|$ and establish conforming decompositions of the overall space,
  thereby satisfying
  \begin{equation}
    V=\widetilde{U}_H \oplus W \quad \text { and } \quad V=\widetilde{V}_H \oplus W.
  \end{equation}
  Furthermore, we have the ‘orthogonality’ relations
  \begin{equation}
    a\left(\widetilde{U}_H, W\right)=0 \quad \text { and } \quad a\left(W, \widetilde{V}_H\right)=0.
  \end{equation}
\end{mylem}
\begin{mypf}
  The proof of this lemma can be found in \cite{RPD2021NHBSS}.
\end{mypf}

\begin{mylem} \label{lemma_inv_ope}
  Let $\mathcal{L}^*$ be the conjugate operator of $\mathcal{L}$, an alternative characterization of the trial and test spaces is provided by
  \begin{equation} \label{eq_LOD_Linv}
    \widetilde{U}_H=\operatorname{span}\left\{\mathcal{L}^{-1} q_j \mid j \in J\right\}
    =\operatorname{span}\left\{\mathcal{L}^{-1} l_j \mid j \in J\right\}
  \end{equation}
  and
  \begin{equation} \label{eq_LOD_Lcinv}
    \widetilde{V}_H=\operatorname{span}\left\{\mathcal{L^*}^{-1} q_j \mid j \in J\right\}
    =\operatorname{span}\left\{\mathcal{L^*}^{-1} l_j \mid j \in J\right\}.
  \end{equation}
\end{mylem}
\begin{mypf}
  The proof of equation (\ref{eq_LOD_Linv}) can be referenced in \cite{RPD2021NHBSS}, and the proof of
  equation (\ref{eq_LOD_Lcinv}) follows a similar approach.
  Noticed that $a(u,v) = \left(\mathcal{L}u, v \right)_{L^2(\Omega)} =\left(u, \mathcal{L^*} v \right)_{L^2(\Omega)}$.
  Let $v = \mathcal{L^*}^{-1} q_j$, then we have $a(w,v) = \left(w, \mathcal{L^*} v \right)_{L^2(\Omega)} = q_j(w) = 0$
  for all $w \in W$, which implies that
  $\operatorname{span}\left\{\mathcal{L^*}^{-1} q_j \mid j \in J\right\} \subseteq \widetilde{V}_H$.
Notice that both spaces possess dimension N, which leads to the conclusion of the lemma.
\end{mypf}

Given the trial and test spaces, the discrete variational problem can be expressed as follows:
Find $\tilde{u}_H \in \widetilde{U}_H$ such that
\begin{equation}
  a\left(\tilde{u}_H, \tilde{v}_H\right)=\left( f,\tilde{v}_H\right)_{L^2(\Omega)} \quad \text { for all }
  \tilde{v}_H \in \widetilde{V}_H.
\end{equation}
\begin{mylem}
  The bases of $\widetilde{U}_H$ and $\widetilde{V}_H$ can be obtained by two set of saddle point problems. For $j \in J$,
  seek $\tilde{u}_k \in V$ and $\mu \in \mathbb{C}^N$, such that
  \begin{equation}
    \begin{aligned}
    a\left(\tilde{u}_k, v\right)+\sum_{j \in J} \mu_j q_j(v) & =0 & & \text { for all } v \in V,\\
    q_j\left(\tilde{u}_k\right) & =\delta_{j k} & & \text { for all } j \in J.
    \end{aligned}
  \end{equation}
  And for $j \in J$, seek $\tilde{v}_k \in V$ and $\nu  \in \mathbb{C}^N$, such that
  \begin{equation}
    \begin{aligned}
    a\left(v,\tilde{v}_k\right)+\sum_{j \in J} \nu_j q_j(v) & =0 & & \text { for all } v \in V,\\
    q_j\left(\tilde{v}_k\right) & =\delta_{j k} & & \text { for all } j \in J.
    \end{aligned}
  \end{equation}
  Then the spaces $\widetilde{U}_H$ and $\widetilde{V}_H$ can be represented as follows
  \begin{equation}
    \widetilde{U}_H = \text{span} \left\{ \tilde{u}_k \right\} \text{ and } \widetilde{V}_H = \text{span}
    \left\{ \tilde{v}_k \right\}.
  \end{equation}
\end{mylem}
\begin{mypf}
  The proof of this lemma can be found in \cite{RPD2021NHBSS}.
\end{mypf}

Indeed, these saddle point problems can be equivalently formulated as energy minimization problems subject to certain constraints. Exploiting this equivalence, the CEM-GMsFEM (Constrained Energy Minimization Generalized Multiscale Finite Element Method) \cite{EYW_CEMGMs_2018} is introduced.

The aforementioned process, commonly referred to as orthogonal decomposition, allows us to construct the trial and test spaces $\widetilde{U}_H$ and $\widetilde{V}_H$ from a general basis function space $U_H$, which is capable of accurately solving the original multiscale problem.
If the basis functions $l_j$ exhibit localization, which means $\text{supp}(l_j)$ represents a small part of the domain $\Omega$, it can be demonstrated that the operators $\mathcal{C}$ and $\mathcal{C}^*$ possess exponential decay properties.
This allows us to localize the projection operators and multiscale basis functions, which is referred to as localized orthogonal decomposition.

Before discussing localization, an interesting question arises: What results can be obtained by repeatedly applying orthogonal decomposition?

\subsection{Iterative orthogonal decomposition}
An iterative sequence of spaces can be constructed based on multiple iterations of orthogonal decomposition, referred to as iterative orthogonal decomposition.
Firstly, we initialize the basis functions
\begin{equation}
u_j^0 := l_j \text{ and } v_j^0 := l_j\text{  for all } j \in J.
\end{equation}
Define
\begin{equation}
  \ U_H^0 := \text{span} \left\{ u_j^0 \right\} \text{ and } V_H^0 := \text{span}
  \left\{ v_j^0 \right\}.
\end{equation}
Then the quantities of interest are defined by
\begin{equation}
  q_j^n := \left(u_j^n, \bullet\right)_{L^2(\Omega)} \text{ and } p_j^n := \left(v_j^n, \bullet\right)_{L^2(\Omega)}
  \text{ for } j \in J, n = 0,1,2,\cdots.
\end{equation}
The kernel spaces are defined by
\begin{equation}
  \begin{aligned}
    W^n:=\left\{v \in V \mid q_j^n(v)=0 \text { for all } j \in J\right\}=\bigcap_{j \in J} \operatorname{ker} q_j^n
    \text{ for } n =0,1,2,\cdots,\\
    X^n:=\left\{v \in V \mid p_j^n(v)=0 \text { for all } j \in J\right\}=\bigcap_{j \in J} \operatorname{ker} p_j^n
    \text{ for } n = 0,1,2,\cdots.
  \end{aligned}
\end{equation}
The projections $\mathcal{C}^n: {V} \rightarrow {W^n}$ and ${\mathcal{C}^*}^n: {V} \rightarrow {W^n}$ are defined by
\begin{equation}
  a(\mathcal{C}^n v, w)=a(v, w) \text{ and } a(w,{\mathcal{C}^*}^n v)=a(w, v) \ \ \forall \ v \in V, w \in W^n,
\end{equation}
for $n = 0,1,2,\cdots$. The trial and test spaces are constructed by
\begin{equation}
  {U}_H^{n+1}:=(1-\mathcal{C}^n) V \quad \text { and } \quad {V}_H^{n+1}:=\left(1-{\mathcal{C}^*}^n\right) V
  \text{ for } n = 0,1,2,\cdots.
\end{equation}

Similarly, the bases of these two sequences of spaces can be constructed by solving two sets of saddle point problems.
For $k \in J, n = 0,1,2,\cdots $,
seek $u_k^{n+1} \in V$ and $\mu^{n+1} \in \mathbb{C}^N$, such that
\begin{equation} \label{eq_saddle_uk}
  \begin{aligned}
    a\left(u_k^{n+1}, v\right)+\sum_{j \in J} \mu_j^{n+1} q_j^n(v) & =0 & & \text { for all } v \in V,\\
    q_j^n \left(u_k^{n+1}\right) & =\delta_{j k} & & \text { for all } j \in J.
  \end{aligned}
\end{equation}
And for $k \in J, n = 0,1,2,\cdots $, seek $v_k^{n+1} \in V$ and $\nu^n  \in \mathbb{C}^N$, such that
\begin{equation}
  \begin{aligned}
    a\left(v,v_k^{n+1} \right)+\sum_{j \in J} \nu_j^{n+1} p_j^n(v) & =0 & & \text { for all } v \in V,\\
    p_j^n\left(v_k^{n+1}\right) & =\delta_{j k} & & \text { for all } j \in J.
  \end{aligned}
\end{equation}
Through iterative orthogonal decomposition, we obtain two sequences of spaces
\begin{equation}
  U_H^0 \rightarrow U_H^1 \rightarrow U_H^2 \rightarrow \cdots \text{ and }  V_H^0 \rightarrow V_H^1 \rightarrow V_H^2
  \rightarrow \cdots.
\end{equation}
Naturally, we can select two spaces from the two sequences for a specific step to serve as the trial and test spaces. However, this method is ineffective as a finite element method for two reasons.
First, the iterative process usually requires a lot of computational resources; second, the obtained functions lack the localizable property.
A very important question is whether limits exist for these two space sequences. If limits do exist, what are their limits?

By Lemma \ref{lemma_inv_ope}, we can deduce
\begin{equation} \label{eq_lemma2_corollary}
  U_H^{n+1} = \operatorname{span} \left\{\mathcal{L}^{-1} u_j^n \mid j \in J\right\} = \mathcal{L}^{-1} U_H^{n}
\end{equation}
and
\begin{equation}
  V_H^{n+1} = \operatorname{span} \left\{\mathcal{L^*}^{-1} v_j^n \mid j \in J\right\} = \mathcal{L^*}^{-1} V_H^{n}
\end{equation}
for $n = 0,1,2,\cdots.$ Then we have
\begin{equation}
  U_H^{n} = \mathcal{L}^{-n} U_H^{0} \text{ and } V_H^{n} = \mathcal{L^*}^{-n} V_H^{0}.
\end{equation}
The above formulas are naturally associated with the power and subspace iteration methods, which are commonly used to solve spectral problems.
Before delving into the discussion of the spectral problem, it is necessary to make some assumptions:
\begin{myassum}
  The linear partial differential operator $\mathcal{L}$ is self-adjoint and positive definite, which means $\mathcal{L} = \mathcal{L^*}$.
\end{myassum}
\begin{myassum} \label{Assum_wellposed}
  The problem (\ref{eq_govering}) is well-posed, indicating the existence of a unique solution $u \in V$ corresponding
  to any given $f \in L^2(\Omega)$.
\end{myassum}
\begin{myassum} \label{Assum_compactemb}
  The Sobolev space $V$ is compactly embedded in $L^2(\Omega)$.
\end{myassum}
Assumption \ref{Assum_wellposed} states that the operator $\mathcal{L}^{-1}$ is a bounded operator from $L^2(\Omega)$ to $V$. Together with Assumption \ref{Assum_compactemb}, it can be proven that the operator $\mathcal{L}^{-1}$ is a compact operator defined on $L^2(\Omega)$. These three assumptions make sure that the operator $\mathcal{L}^{-1}$ is both compact and self-adjoint. This allows us to use the spectral decomposition theorem, which is designed to work with compact self-adjoint operators \cite{Brezis_FASSPDE_2011}.

Let $\left\{ \left( \lambda_i,\phi_i\right) | i = 0,1,2,\cdots \right\}$
denote the set of all eigenpairs of the inverse operator $\mathcal{L}^{-1}$ in $V$, where
\begin{equation}
  \lambda_0 \geq \lambda_1 \geq \lambda_2 \geq \cdots \geq 0.
\end{equation}
Consequently, it can be deduced that the set $\left\{ \phi_i \right\}$ forms a set of complete orthogonal bases of $V$.
We can make the assumption that $u_j^0 \in V \text{ for all } j \in J$.
In the event that this assumption is not met, we redefine $u_j^0 := u_j^1 \in V$, ensuring that the iterative sequence remains unaltered.
\begin{mythm} \label{theorem_limit}
  Let $P_i$ is the spectral projector associated with the eigenvalues $\lambda_1,\lambda_2,\cdots,\lambda_i$. Assume that
  $\operatorname{rank}\left(P_i\left[u_1^0, u_2^0, \ldots, u_i^0\right]\right)=i$ for $i=1,2, \ldots, N$ and
  $\lambda_N > \lambda_{N+1}$. We can conclude that the sequence of spaces $\left\{ U_H^n\right\}$ converges to the
  eigensubspace $V_{\text{eig}}$ spanned by the first $N$ eigenfunctions of the operator $\mathcal{L}^{-1}$.
\end{mythm}
\begin{mypf}
  Let $u_j^0 = \sum_{i=1}^{\infty} \zeta_j^i \phi_i $ for $j \in J$. Based on the assumption, there exists another set of
  functions $\left\{ \bar{u}_j^0\right\}$ that can be expressed as linear combinations of $\left\{ {u}_j^0\right\}$, satisfying
  \begin{equation}
    \bar{u}_j^0 = \sum_{i=j}^{\infty} \bar{\zeta}_j^i \phi_i, \text{ for } j \in J,
  \end{equation}
  where $\bar{\zeta}_j^j \neq 0$. Then we have
  \begin{equation}
    \mathcal{L}^{-n} \bar{u}_j^0 = \sum_{i=j}^{\infty} \bar{\zeta}_j^i \lambda_i^n \phi_i, \text{ for } j \in J.
  \end{equation}
  After normalization,
  \begin{equation}
    \hat{u}_j^n :=\frac{1}{\lambda_j^n}\mathcal{L}^{-n} \bar{u}_j^0 = \bar{\zeta}_j^j \phi_j + \sum_{i=j+1}^{\infty} \bar{\zeta}_j^i
    \left(\frac{\lambda_i}{\lambda_j} \right)^n \phi_i, \text{ for } j \in J.
  \end{equation}
  It is easy to verify that $U_H^n = \operatorname{span} \left\{ \hat{u}_j^n \right\}$. If $\lambda_j >\lambda_{j+1}$, we have
  $\lim_{n \to \infty} \hat{u}_j^n =  \bar{\zeta}_j^j \phi_j$. If $\lambda_j =\lambda_{j+1}= \cdots =\lambda_{j+k} > \lambda_{j+k+1}$, we have
  $\lim_{n \to \infty} \hat{u}_{j+l}^n =  \sum_{i=j+l}^{j+k}\bar{\zeta}_{j+l}^i \phi_i$ for $l = 0,1,\cdots,k$. In conclusion,
  \begin{equation}
    \lim_{n \to \infty} U_H^n = \operatorname{span} \left\{\lim_{n \to \infty} \hat{u}_j^n | j \in J\right\} =
    \operatorname{span} \left\{ \phi_j | j \in J \right\}.
  \end{equation}
  the proof is completed.
\end{mypf}

In fact, orthogonal decomposition can be viewed as a special case of the subspace iteration method.
Start with $U_H^n = \left\{ u_j^n\right\}$, apply the operator $\mathcal{L}^{-1}$,  and we obtain $\left\{ \mathcal{L}^{-1} u_j^n\right\}$.
In traditional subspace iteration methods, some orthogonalization techniques are used to obtain new basis functions, such as the Gram-Schmidt orthogonalization method. In orthogonal decomposition, based on the saddle problem (\ref{eq_saddle_uk}), we seek
\begin{equation}
  u_k^{n+1} = \sum_{j \in J} -\mu_j^{n+1} \mathcal{L}^{-1} u_j^n,
\end{equation}
such that $\left( u_k^{n+1},u_j^{n}\right)_{L^2(\Omega)} = \delta_{jk}$.
This can be considered as a special orthogonalization method, different from the traditional orthogonalization methods that satisfy $\left( u_k^{n+1},u_j^{n+1}\right)_{L^2(\Omega)} = \delta_{jk}$.

\begin{myrema}
  Orthogonalization techniques in the traditional subspace iteration methods are executed sequentially, facilitating the convergence of the basis functions towards their respective eigenfunctions. In the iterative orthogonal decomposition, the sequence of functions $\left\{u_j^n \right\}_{n=1}^{\infty}$ does not converge.
  It is worth noting that the sequence of functions can be decomposed into two alternating subsequences, each of which exhibits convergence, based on empirical observations from numerical experiments.
\end{myrema}

\section{Two subspaces of spectral problem algorithms} \label{section:num3}
In the preceding section, we provide a novel perspective that some multiscale methods can be be regarded as one iteration step under approximating the eigenspace of the corresponding local spectral problems.
In fact, a multitude of multiscale methodologies exhibit a strong interconnection with spectral problems \cite{YJH_GMsFEM_2013,IR_OLASfGFEM_2011,EYW_CEMGMs_2018}. This impetus drives us to delve into broader possibilities within multiscale modeling, commencing with  spectral problem algorithms.

We consider a spectral problem
\begin{equation} \label{eq_lineareig}
  A x = \lambda x,
\end{equation}
where $A \in \mathbb{C}^{m \times m}$ is hermitian. Let us assume that the eigenvalues are arranged in descending order, which means
\begin{equation}
  \lambda_1 \geq \lambda_2 \geq \cdots \geq \lambda_m.
\end{equation}
A classic mission is to find the leading $p$ eigenpairs $\left(\lambda_i,x_i\right)$ for $i=1,2,\cdots,p$. To achieve this, a natural choice is the subspace iteration method, a foundational and simple method.
\subsection{Standard subspace iteration}
Algorithm \ref{alg:subiter} shows the standard subspace iteration, which is capable of approximating the leading $p$ eigenpairs $\left(\lambda_i,x_i\right), i=1,2,\cdots,p$. The convergence speed for each corresponding eigenfunction $x_i$ can be expressed as $\left(\left|\frac{\lambda_{p+1}}{\lambda_i}\right|+\epsilon_n\right)^n$, where $n$ denotes the number of iteration steps and $\epsilon_n$ tends to zero \cite{SY_NMfLEP_2011}.
\begin{algorithm}[H]
   \caption{Standard Subspace Iteration}
  \label{alg:subiter}
  \begin{algorithmic}
    \STATE 1. \textbf{Start:} Select an initial set of vectors $X_0=\left[x_1^0, \ldots, x_p^0\right]$.\\
    \STATE 2. \textbf{Iterate:} Repeat until convergence is achieved,\\
    \ (a) Calculate $X_k:=A X_{k-1}$\\
    \ (b) Calculate $X_k=Q R$ the $Q R$ factorization of $X_k$, and set $X_k:=Q$.
 \end{algorithmic}
\end{algorithm}
The standard subspace iteration is the most fundamental algorithm for solving spectral problems. Compared to other methods, it may not offer a significant advantage. However, its significance lies in its ability to provide insights into multiscale method formulation, which we will discuss in the next section. Furthermore, applying specific projection or preprocessing strategies has the potential to accelerate computational processes.

\subsection{Krylov subspace iteration}
The Krylov subspace iteration, as a trivial extension of the standard subspace iteration, is one of the most important methods available for computing the eigenvalues and eigenvectors of large matrices, particularly in the Hermitian case.
In comparison to the standard subspace iteration, the Krylov subspace iteration demonstrates enhanced efficiency, memory utilization, and flexibility.

For the same spectral problem (\ref{eq_lineareig}), the Krylov subspace is defined by
\begin{equation}
  \mathcal{K}_r(A,x) := \operatorname{span}\left\{x, A x, A^2 x, \ldots A^{r-1} x\right\}.
\end{equation}
If there is no possibility of ambiguity, $\mathcal{K}_r (A,x)$ is denoted as $\mathcal{K}_r$. In contrast with the standard subspace iteration, the Krylov subspace iteration necessitates only a single matrix or operator operation at each iteration step, as opposed to multiple operations.
A few well-known of these Krylov subspace methods are Arnoldi's method and Lanczos' method.
This paper will introduce Arnoldi's method as an illustrative example.

Arnoldi's method is an orthogonal projection method onto $\mathcal{K}_r$ for large sparse matrices. 
Specifically, the procedure introduced by Arnoldi in 1951 starts by building an orthogonal basis of the Krylov subspace $\mathcal{K}_r$.
Subsequently, we approximate eigenpairs within the subspace $\mathcal{K}_r$ using orthogonal projection techniques, such as the Rayleigh-Ritz procedure.
There are several distinct implementations of Arnoldi’s method, which are all mathematically equivalent, and Algorithm \ref{alg:ArnoldiKrylov} is one of them.
\begin{algorithm}[H]
  \caption{Arnoldi's Krylov subspace iteration}
 \label{alg:ArnoldiKrylov}
 \begin{algorithmic}
   \STATE 1. \textbf{Start:} Select an initial vector $x_1$.\\
   \STATE 2. \textbf{Iterate:} Compute for $j = 1,2,\cdots,l-1$:\\
   \ (a) $h_{i j}=\left(A x_j, x_i\right), \quad i=1,2, \ldots, j$, \\
   \ (b) $y_j=A x_j-\sum_{i=1}^j h_{i j} x_i$, \\
   \ (c) $h_{j+1, j}=\left\|y_j\right\|_2, \quad \text { if } h_{j+1, j}=0 \text { stop }$ \\
   \ (d) $x_{j+1}=y_j / h_{j+1, j}$.
   \STATE 3. \textbf{Orthogonal projection:} Use the Rayleigh-Ritz procedure to obtain desired eigenpairs in $\mathcal{K}_r = \operatorname{span} \left\{x_1,x_2,\cdots,x_l\right\}$.
\end{algorithmic}
\end{algorithm}
Indeed, there are numerous efficient methods for spectral problems, such as the Jacobi-Davidson method \cite{Sleijpen_JacobiDavidson_2000}, and this article only highlights a few of the most common methods.

\section{From spectral problem to LSI} \label{section:num4}
Virtually every iteration method used for spectral problems has the potential to be extended to a multiscale modeling method.
Specifically, we can extract a step or multiple steps from the iterative process to construct a subspace that approximates the eigenfunction subspace.
Based on the standard subspace iteration method, we proposed the localized standard subspace iteration (LSSI) method. Furthermore, based on the Krylov subspace iteration, we proposed the localized Krylov subspace iteration (LKSI) method. Before introducing these methods, we start with a localization process, analogous to the technique adopted in the majority of multiscale methods.


Let $\left\{\omega_i \right\}_{i=1}^{N_c}$ be a open cover of domain $\Omega$, such that $\Omega = \bigcup_{i=1}^{N_c}  \omega_i$, where $N_c$ is the number of subdomains. An elementary choice for the set $\left\{\omega_i \right\}$ is to extend each element $K_i$ of the finite element partition $\mathcal{T}_H$ by one or several layers. There is a set of partition of unity $\left\{\chi_i \right\}_{i=1}^{N_c}$, such that
\begin{equation}
  1 = \sum_{i=1}^{N_c} \chi_i \text{ and } \operatorname{supp}(\chi_i) =  \omega_i \text{ for } i =1,2,\cdots,N_c.
\end{equation}
Let $V(\omega_i) := \{ v \in V | \operatorname{supp}(v) \in \omega_i \}$. The operator $\mathcal{L}$ restricted to the space $V(\omega_i)$ is denoted as $\mathcal{L}_i$.
It is clear that $\mathcal{L}_i^{-1}$, being defined on $L^2(\omega_i)$, is also compact and self-adjoint. Therefore, we can define local spectral problem
\begin{equation} \label{eq:local:spectral}
  \mathcal{L}_i^{-1} \phi_i^{j} = \lambda_i^j \phi_i^{j},   \quad \phi_i^{j} \in V\left(\omega_i\right),
\end{equation}
for $i = 1,2,\cdots,N_c, j = 1,2, \cdots $.

\subsection{Localized standard subspace iteration (LSSI) method}
For each subdomain $\omega_i$, there is a set of elementary basis functions denoted as $\phi_i^{j,0} \in V\left(\omega_i\right)$, where $j = 1,2,\cdots,L_i$ and $L_i$ represents the number of basis functions in subdomain $\omega_i$.
Through the standard subspace iteration in each subdomain, a sequence of basis function sets can be acquired. To be specific, for $i = 1,2,\cdots,N_c, k = 1,2,\cdots,L_i, n = 0,1,2,\cdots$, we seek $\phi_i^{k,n+1} \in V\left(\omega_i\right)$ and $\mu_i^{n+1} \in \mathbb{C}^{L_i}$, such that
\begin{equation} \label{eq_saddle_SMs}
  \begin{aligned}
    a\left(\phi_i^{k,n+1}, v\right)+\sum_{j =1}^{L_i} \mu_i^{j,n+1} q_i^{j,n}(v) & =0 & & \text { for all } v \in V\left(\omega_i\right),\\
    q_i^{j,n} \left(\phi_i^{k,n+1}\right) & =\delta_{j k} & & \text { for all } j =1,2,\cdots,L_i.
  \end{aligned}
\end{equation}
In the above equation, $ q_i^{j,n} $ is a linear functional defined by
\begin{equation}
  q_i^{j,n} := \left(\phi_i^{j,n}, \bullet\right)_{L^2(\Omega)}.
\end{equation}
The local multiscale space $V_{\text{S}}^{i,n}$ is constructed by
\begin{equation} \label{eq:local:multi:space}
  V_{\text{S}}^{i,n} := \text{span} \{\phi_i^{j,n} | j = 1,2,\cdots L_i\},
\end{equation}
$\text{ for } i = 1,2,\cdots,N_c, n = 0,1,2,\cdots$. The multiscale space $V_{\text{S}}^{n}$ is constructed by
\begin{equation}
  V_{\text{S}}^{n} := \bigoplus_{i = 1}^{N_c} V_{\text{S}}^{i,n}.
\end{equation}
Similar to Eq.(\ref{eq_lemma2_corollary}), it can be deduced that
\begin{equation} \label{eq_lemma2_corollary_local}
  V_{\text{S}}^{i,n+1} = \operatorname{span} \left\{\mathcal{L}_i^{-1} \phi_i^{j,n} \mid j =1,2,\cdots,L_i\right\} = \mathcal{L}_i^{-1} V_{\text{S}}^{i,n} = \mathcal{L}_i^{-(n+1)} V_{\text{S}}^{i,0}.
\end{equation}
As a direct corollary of Theorem \ref{theorem_limit},
\begin{equation}
  \lim_{n \to \infty} V_{\text{S}}^{i,n} = \operatorname{span} \left\{ \phi_i^j | j = 1,2,\cdots,L_i\right\}.
\end{equation}
Once the multiscale space $V_{\text{S}}^{n}$ is obtained, the multiscale solution can be obtained using the Galerkin method. Due to its strong resemblance to the standard subspace iteration in spectral problem algorithms, this method is referred to as the localized standard subspace iteration (LSSI) method. The convergence rate of the LSSI depends on the separation of eigenvalues. Given the rapid decay characteristic of eigenvalues in multiscale spectral problems, it often only takes a few iterations to achieve highly satisfactory results.
\subsection{Localized Krylov subspace iteration (LKSI) method}
When employing the subspace iteration method to solve eigenvalue problems, it is common practice to exclusively use the result of the final iteration step for computation. In reality, the functions obtained at each iteration step can all be used for eigenfunction calculation. This fundamental principle reveals the core of the Krylov method, wherein it exhibits a distinct advantage in terms of computational efficiency and memory usage.

In each subdomain $\omega_i$, suppose there is an initial basis function $\psi_i^{0} \in V\left(\omega_i\right)$. For $i = 1,2,\cdots,N_c, n = 0,1,\cdots$, we seek $\psi_i^{n+1}$ and $\mu_i^{n+1} \in \mathbb{C}$ , such that
\begin{equation} \label{eq_saddle_KSMs}
  \begin{aligned}
    a\left(\psi_i^{n+1}, v\right)+ \mu_i^{n+1} q_i^{n}(v) & =0 & & \text { for all } v \in V\left(\omega_i\right),\\
    q_i^{n} \left(\psi_i^{n+1}\right) & =1.
  \end{aligned}
\end{equation}
Similarly, the definition of $q_i^{n}$ is as follows
\begin{equation}
  q_i^{n} := \left(\psi_i^{n}, \bullet\right)_{L^2(\Omega)}.
\end{equation}
The local multiscale space $V_{\text{K}}^{i,n}$ is constructed by
\begin{equation} \label{eq:ksms:local:space}
  V_{\text{K}}^{i,n} := \operatorname{span} \{\psi_i^{k} | k = 0,1,2,\cdots n \},
\end{equation}
$\text{ for } i = 1,2,\cdots,N_c, n = 0,1,2,\cdots$. The multiscale space $V_{\text{K}}^{n}$ is constructed by
\begin{equation}
  V_{\text{K}}^{n} := \bigoplus_{i = 1}^{N_c} V_{\text{K}}^{i,n}.
\end{equation}
It is easy to deduce that
\begin{equation}
  V_{\text{K}}^{i,n} = \operatorname{span} \left\{\mathcal{L}_i^{-k} \psi_i^{0} \mid k =0,1,2,\cdots,n \right\}.
\end{equation}
In fact, $V_{\text{K}}^{i,n} $ is the Krylov subspace of the operator $\mathcal{L}_i^{-1}$ with the initial function $\psi_i^{0}$.
\section{Convergence analysis}  \label{section:num5}
In this section, we will establish the convergence of our proposed methods.
In order to demonstrate the convergence of the proposed methods, we will initially show the interpolation error of using the local eigenfunctions as basis functions.

\subsection{Interpolation error}
Define the local eigenfunction space $V_{\text{eig}}$ as follows:
\begin{equation}
  V_{\text{eig}} := \operatorname{span} \{ \phi_i^j | i = 1,2,\cdots,N_c, j =  1,2,\cdots,L_i\},
\end{equation}
where $\phi_i^j$ is the local eigenfunction difined by Eq. \cref{eq:local:spectral}.
Then for all $u \in V$, define interpolation operator $\mathcal{I}_{\text{eig}}: V \rightarrow V_{\text{eig}}$ as follows:
\begin{equation}
  \mathcal{I}_{\text{eig}} u := \sum_{i=1}^{N_c} \sum_{j=1}^{L_i} \left \langle \chi_i u, \phi_i^j \right \rangle \phi_i^j.
\end{equation}

\begin{mythm} \label{theorem_interp_spec}
  Suppose $u \in V$ , then we have an estimation for the interpolation error
  \begin{equation} \label{eq:theorem1:1}
    \Vert u - \mathcal{I}_{\text{eig}} u \Vert_a \leq \sqrt{\lambda^{L+1}}  \sum_{i=1}^{N_c} \Vert \mathcal{L}\chi_i u \Vert_{L^2(\Omega)},
  \end{equation}
  where $\lambda^{L+1} := \mathop{\max} \limits_{i} \ \lambda_i^{L_i+1}$, and the energy norm $\Vert \bullet \Vert_a$ is defined by $\Vert u \Vert_a^2 = a(u,u)$.
\end{mythm}
\begin{mypf}
  Noticed that
  \begin{equation}
    \begin{aligned}
    u - \mathcal{I}_{\text{eig}} u & = \sum_{i=1}^{N_c} \chi_i u - \sum_{i=1}^{N_c} \sum_{j=1}^{L_i} \left \langle \chi_i u, \phi_i^j \right \rangle \phi_i^j \\
    & = \sum_{i=1}^{N_c} \sum_{j=L_i+1}^{\infty} \left \langle \chi_i u, \phi_i^j \right \rangle \phi_i^j .
  \end{aligned}
  \end{equation}
  It is easy to verify that
  \begin{equation}
    \begin{aligned}
      \Vert \sum_{j=L_i+1}^{\infty} \left \langle \chi_i u, \phi_i^j \right \rangle \phi_i^j \Vert_a^2 &  = \sum_{j=L_i+1}^{\infty} \left \langle \chi_i u, \phi_i^j \right    \rangle^2 \Vert \phi_i^j \Vert_a^2    \\
      & =    \sum_{j=L_i+1}^{\infty} \left \langle \chi_i u, \phi_i^j \right \rangle^2  \frac{1}{\lambda_i^j} \\
      & \leq \lambda_i^{L_i+1} \sum_{j=L_i+1}^{\infty} \left \langle \chi_i u, \phi_i^j \right \rangle^2  \left(\frac{1}{\lambda_i^j}\right)^2\\
      & \leq \lambda_i^{L_i+1} \sum_{j=1}^{\infty} \left \langle \chi_i u, \phi_i^j \right \rangle^2  \left(\frac{1}{\lambda_i^j}\right)^2\\
      & = \lambda_i^{L_i+1} \Vert \mathcal{L}\chi_i u \Vert_{L^2(\Omega)}^2 .
    \end{aligned}
  \end{equation}
The last equation is based on the expansion of the operator $\mathcal{L}$ on $V(\omega_i)$,
\begin{equation}
  \mathcal{L} \chi_i u = \sum_{j=1}^{\infty} \frac{1}{\lambda_i^j} \left \langle \chi_i u, \phi_i^j \right \rangle \phi_i^j.
\end{equation}
Then we have
\begin{equation}
\begin{aligned}
  \Vert u - \mathcal{I}_{\text{eig}} u \Vert_a & \leq \sum_{i=1}^{N_c} \Vert \sum_{j=L_i+1}^{\infty} \left \langle \chi_i u, \phi_i^j \right \rangle \phi_i^j \Vert_a \\
  & \leq \sum_{i=1}^{N_c} \sqrt{\lambda_i^{L_i+1}} \Vert \mathcal{L}\chi_i u \Vert_{L^2(\Omega)}  \\
  & \leq \sqrt{\lambda^{L+1}}  \sum_{i=1}^{N_c} \Vert \mathcal{L}\chi_i u \Vert_{L^2(\Omega)} .
\end{aligned}
\end{equation}

\end{mypf}

\subsection{Convergence result of the LSSI method}
Before the convergence result of the LSSI is obtained, let us give a lemma without proof.
\begin{mylem} \label{lemma:first}
  Suppose $\phi_i^j$ is the $j$th eigenfunction of the local spectral problem \cref{eq:local:spectral}, and the local multiscale space $V_{\text{S}}^{i,n}$ is defined by \cref{eq:local:multi:space}. There exists a $\hat{\phi}_i^{j,n}$ in $V_{\text{S}}^{i,n}$, such that the following inequality is satisfied:
  \begin{equation}
    \Vert \hat{\phi}_i^{j,n} -  \phi_i^j \Vert_a \leq C(\phi_i^j,V_{\text{S}}^{i,0}) \left( \frac{\lambda_i^{L_i+1}}{\lambda_i^{j}}+ \epsilon_n \right)^n,
  \end{equation}
  where $\epsilon_n$ tends to zero as n tends to infinity.
\end{mylem}

In fact, Let $P_i$ is the spectral projector associated with the invariant subspace associated with $\lambda_i^1,\dots,\lambda_i^{L^i}$. Then for each $\phi_i^j$, there exists a unique $s \in V_{\text{S}}^{i,0}$ such that $P_i s = \phi_i^j$. The constant $C(\phi_i^j,V_{\text{S}}^{i,0})$ is defined by
\begin{equation}
  C(\phi_i^j,V_{\text{S}}^{i,0})  = \Vert s -\phi_i^j \Vert_a.
\end{equation}
For the sake of simplicity, $C(\phi_i^j,V_{\text{S}}^{i,0})$ is replaced by a constant $C$ in the following. For details of the proof, we refer the reader to Yousef's book (\cite{SY_NMfLEP_2011}; Theorem 5.2).

Based on the auxiliary function $\hat{\phi}_i^{j,n}$ in $V_{\text{S}}^{i,n}$, we define Interpolation operator $\mathcal{I}_{S}: V\rightarrow V_{\text{S}}^n$,
\begin{equation}
  \mathcal{I}_{\text{S}} u := \sum_{i=1}^{N_c} \sum_{j=1}^{L_i} \left \langle \chi_i u, \phi_i^j \right \rangle \hat{\phi}_i^j.
\end{equation}

\begin{mythm}
Suppose $u$ is the solution of \cref{eq_govering}, and $u_{\text{S}}$ is the finite element solution in the multiscale space $V_{\text{S}}^n$. The following convergence result holds:
  \begin{equation} \label{eq:theorem2:1}
    \Vert u- u_{S} \Vert_a \leq \sqrt{\lambda^{L+1}}  \sum_{i=1}^{N_c} \Vert \mathcal{L}\chi_i u \Vert_{L^2(\Omega)} +C \left( \lambda^{\frac{L+1}{L}}+ \epsilon_n \right)^n  \Vert  u \Vert_{L^2(\Omega)},
  \end{equation}
  where $\lambda^{\frac{L+1}{L}}$ is defined by
  \begin{equation}
    \lambda^{\frac{L+1}{L}} := \mathop{\max} \limits_i \frac{\lambda_i^{L_i+1}}{\lambda_i^{L_i}}.
  \end{equation}
\end{mythm}
\begin{mypf}
  By the c\'ea Lemma and triangular inequality, we have
  \begin{equation}
    \begin{aligned}
      \Vert u- u_{S} \Vert_a  & = \mathop{\inf} \limits_{v \in V_{\text{S}}^n} \Vert u- v \Vert_a
       \leq \Vert u- \mathcal{I}_{\text{S}} u \Vert_a \\
      & \leq \Vert u- \mathcal{I}_{\text{eig}} u \Vert_a + \Vert \mathcal{I}_{\text{eig}} u- \mathcal{I}_{\text{S}} u \Vert_a.
    \end{aligned}
  \end{equation}
  Based on the definitions of $\mathcal{I}_{\text{eig}}$ and $\mathcal{I}_{\text{S}}$, we have
  \begin{equation}
    \begin{aligned}
      \mathcal{I}_{\text{eig}} u- \mathcal{I}_{\text{S}} u = \sum_{i=1}^{N_c} \sum_{j=1}^{L_i} \left \langle \chi_i u, \phi_i^j \right \rangle \left( \phi_i^j-\hat{\phi}_i^j \right).
    \end{aligned}
  \end{equation}
  Using the \cref{lemma:first},
  \begin{equation} \label{eq:theorem2:3}
    \begin{aligned}
    \Vert \mathcal{I}_{\text{eig}} u- \mathcal{I}_{\text{S}} u \Vert_a & \leq  \sum_{i=1}^{N_c} \sum_{j=1}^{L_i} \left \langle \chi_i u, \phi_i^j \right \rangle \Vert \hat{\phi}_i^{j,n} -  \phi_i^j \Vert_a \\
    & \leq C \sum_{i=1}^{N_c} \sum_{j=1}^{L_i} \left \langle \chi_i u, \phi_i^j \right \rangle \left( \frac{\lambda_i^{L_i+1}}{\lambda_i^{j}}+ \epsilon_n \right)^n \\
    & \leq C \sum_{i=1}^{N_c}  \left( \frac{\lambda_i^{L_i+1}}{\lambda_i^{L_i}}+ \epsilon_n \right)^n \sum_{j=1}^{L_i} \left \langle \chi_i u, \phi_i^j \right \rangle \\
    & \leq C \left( \lambda^{\frac{L+1}{L}}+ \epsilon_n \right)^n \sum_{i=1}^{N_c} \Vert \chi_i u \Vert_{L^2(\Omega)} \\
    & \leq C \left( \lambda^{\frac{L+1}{L}}+ \epsilon_n \right)^n \Vert  u \Vert_{L^2(\Omega)}
    \end{aligned}
  \end{equation}
  Combining Eqs. \cref{eq:theorem1:1} and \cref{eq:theorem2:3}, the proof is completed.
\end{mypf}

\subsection{Convergence result of the LKSI method}
Similar to \cref{lemma:first}, we first give the following lemma without proof.
\begin{mylem} \label{lemma:second}
  Suppose $\phi_i^j$ is the $j$th eigenfunction of the local spectral problem \cref{eq:local:spectral}, and the local multiscale space $V_{\text{K}}^{i,n}$ is defined by \cref{eq:ksms:local:space}. There exists a $\bar{\phi}_i^{j,n}$ in $V_{\text{K}}^{i,n}$, such that the following inequality is satisfied
  \begin{equation}
    \Vert \bar{\phi}_i^{j,n} -  \phi_i^j \Vert_a \leq C \tan \theta \left( \phi_i^j, \psi_i^0 \right) \frac{\alpha_i^j}{\left(1+4 \gamma_i^j\right)^{n-j}} ,
  \end{equation}
  where
  \begin{equation}
    \alpha_i^1 = 1, \ \alpha_i^j = \prod \limits_{k=1}^{j-1} \frac{\lambda_i^k}{\lambda_i^k-\lambda_i^j} \text{ for } j >1,
  \end{equation}
  and
  \begin{equation}
    \gamma_i^j = \frac{\lambda_i^j-\lambda_i^{j+1}}{\lambda_i^{j+1}}.
  \end{equation}
\end{mylem}
For the sake of simplicity, $C \tan \theta \left( \phi_i^j, \psi_i^0 \right) \frac{\alpha_i^j}{\left(1+4 \gamma_i^j\right)^{-j}} $ is replaced by a constant $C$ in the following. For details of the proof, we also refer the reader to Yousef's book (\cite{SY_NMfLEP_2011}; Theorems 4.8 and 6.3).

Based on the auxiliary function $\bar{\phi}_i^{j,n}$ in $V_{\text{K}}^{i,n}$, we define Interpolation operator $\mathcal{I}_{\text{K}}: V\rightarrow V_{\text{K}}^n$,
\begin{equation}
  \mathcal{I}_{\text{K}} u := \sum_{i=1}^{N_c} \sum_{j=1}^{L_i} \left \langle \chi_i u, \phi_i^j \right \rangle \bar{\phi}_i^j.
\end{equation}

\begin{mythm}
  Suppose $u$ is the accuracy solution of \cref{eq_govering}, $u_{\text{K}}$ is the finite element solution in the multiscale space $V_{\text{K}}^n$. The following convergence result holds:
  \begin{equation}
    \Vert u- u_{K} \Vert_a \leq \sqrt{\lambda^{L+1}}  \sum_{i=1}^{N_c} \Vert \mathcal{L}\chi_i u \Vert_{L^2(\Omega)} +C \left(\frac{1}{1+4\Gamma} \right)^n  \Vert u \Vert_{L^2(\Omega)},
  \end{equation}
  where $\Gamma$ is defined by
  \begin{equation}
    \Gamma := \mathop{\min} \limits_{i,j} \gamma_i^j.
  \end{equation}
\end{mythm}
The proof of this theorem is similar to the previous theorem, so we won’t go into details here.

\section{Numerical experiments}  \label{section:num6}
In this section, we present several numerical examples to evaluate the performance of the proposed LSSI and LKSI methods. The methods proposed in this article are effective for symmetric positive-definite differential operators.

For each element $K_i$ of the coarse finite element mesh $\mathcal{T}_H$, define the oversampling block $K_i^{m} $ as follows,
\begin{equation}
  \begin{aligned}
    & K_i^{0} := K_i,\\
    & K_i^{m} :=  \text{int} \left( \bigcup T \in \mathcal{T}_H|  T\cap \overline{K}_i^{m-1} \neq  \emptyset \right), m = 1,2,3... .
  \end{aligned}
\end{equation}
In this work, the subdomain $\omega_i$ is chosen as $K_i^{m} $, where $m$ is the number of oversampling layers. In addition, a fine finite element mesh $\mathcal{T}_h$ is used to obtain the reference solution and solve local problems. The errors in the following numerical examples are relative errors compared to the reference solution. Our numerical examples are performed on a desktop workstation with 16 GB of memory and a 3.4GHz Core i7 CPU.
\subsection{Diffusion problem} \label{sec_NumExam_ex1}
\begin{figure}[htbp]
  \centering
  \includegraphics[width=0.5\linewidth] {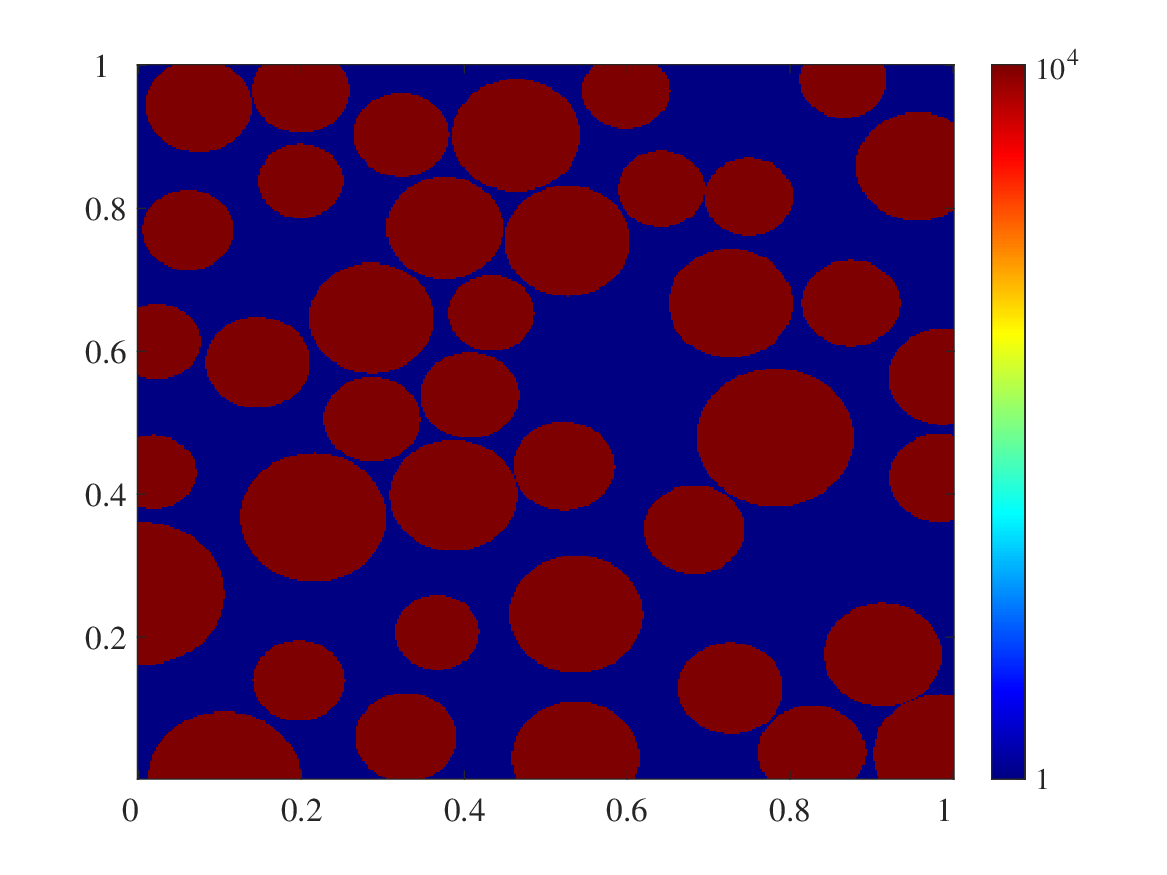}
  \caption{The high-contrast permeability coefficient $\kappa$.}
  \label{ex1_kappa}
\end{figure}

We consider a diffusion problem
\begin{equation} \label{eq_diffusion}
  - \nabla \cdot \left( \kappa \nabla u \right) = f \text{ in } \Omega,
\end{equation}
where $\Omega = [0,1]^2$ and $f = \sin(\pi x) \sin(\pi y) $. $\kappa$ is a high-contrast permeability coefficient with multiscale characteristics, and is shown in ~\cref{ex1_kappa}. The fine mesh size is $h=1/100$, and the coarse mesh size is $H=1/10$. In the LOD and LSSI, we select $4$ bilinear functions on $K_i$ as the initial basis functions $\left\{ \phi_i^{j,0} \right\}_{j=1}^{4}$ for each subdomain $K_i^{m}$. In the LKSI, we select a piecewise constant on $K_i$ as the initial basis function $\psi_i^{0} $ for each subdomain $K_i^{m}$.
We use the notation `$LSSI\text{-}n$' to represent the LSSI and the `$LKSI\text{-}n$' to represent the LKSI with $n$ iteration steps.

\begin{figure}[htbp]
  \centering
  {\tiny(a)}\includegraphics[width=0.45\textwidth]{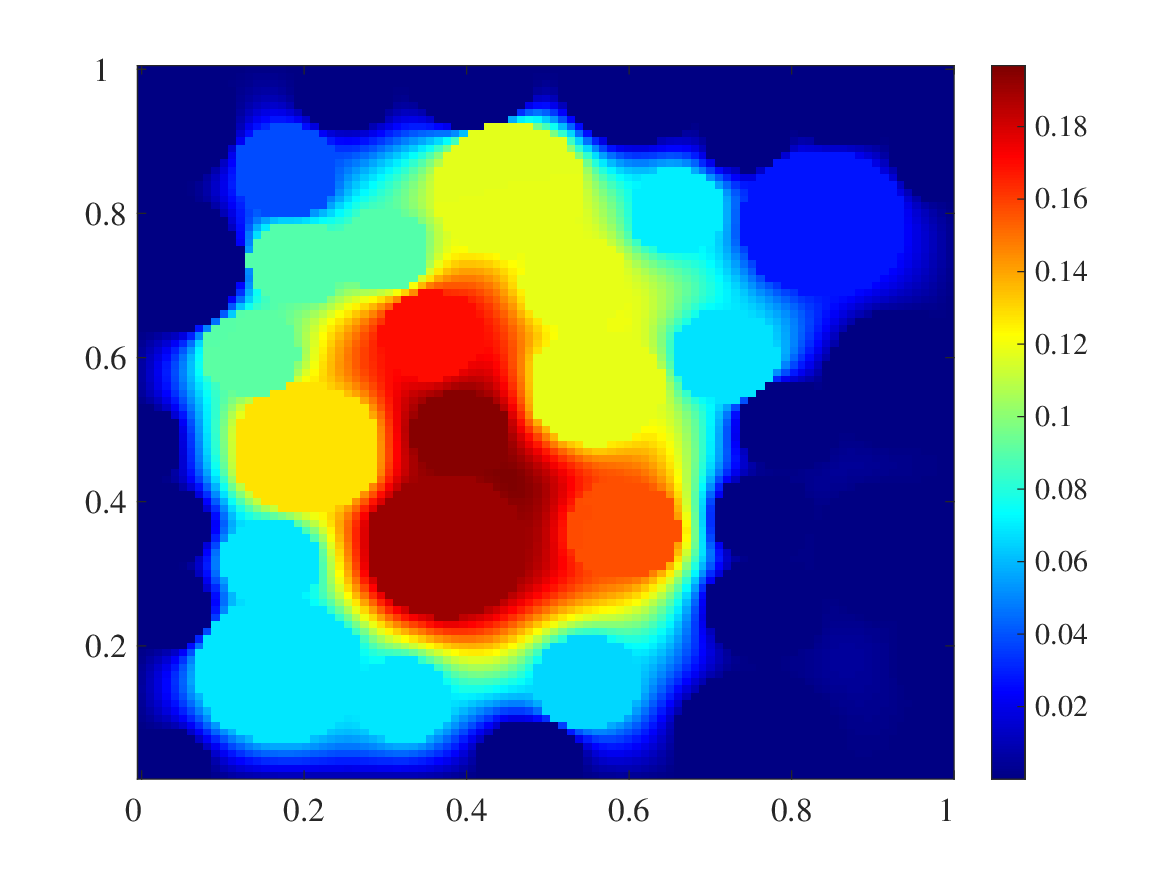}~
  {\tiny(d)}\includegraphics[width=0.45\textwidth]{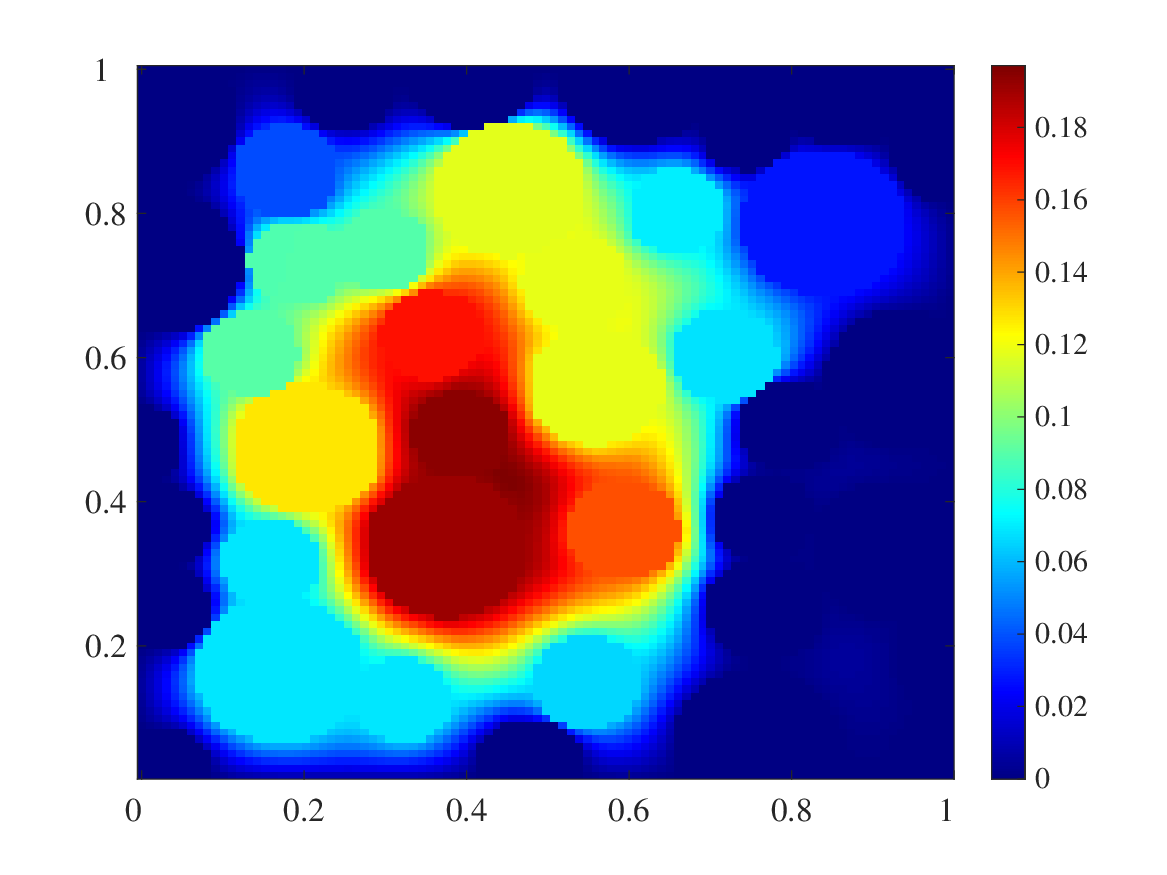}
  {\tiny(b)}\includegraphics[width=0.45\textwidth]{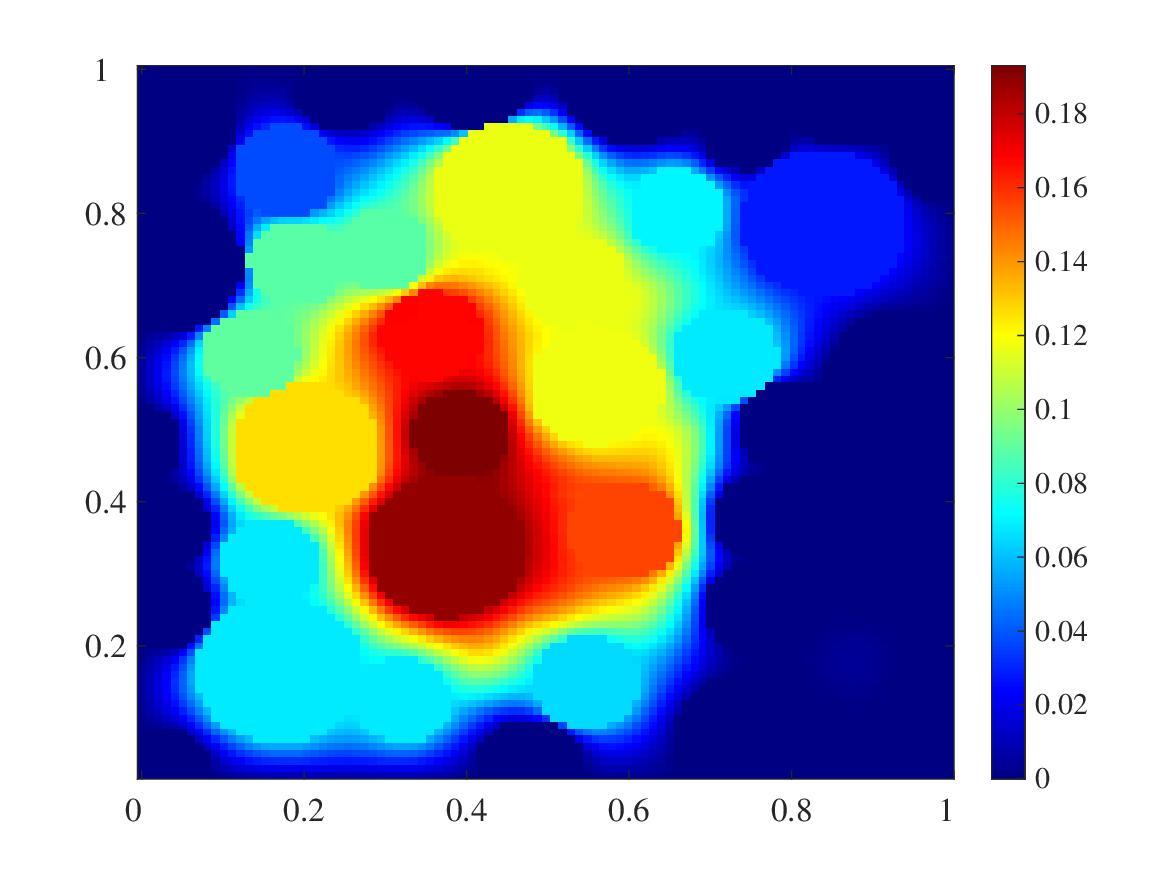}~
  {\tiny(e)}\includegraphics[width=0.45\textwidth]{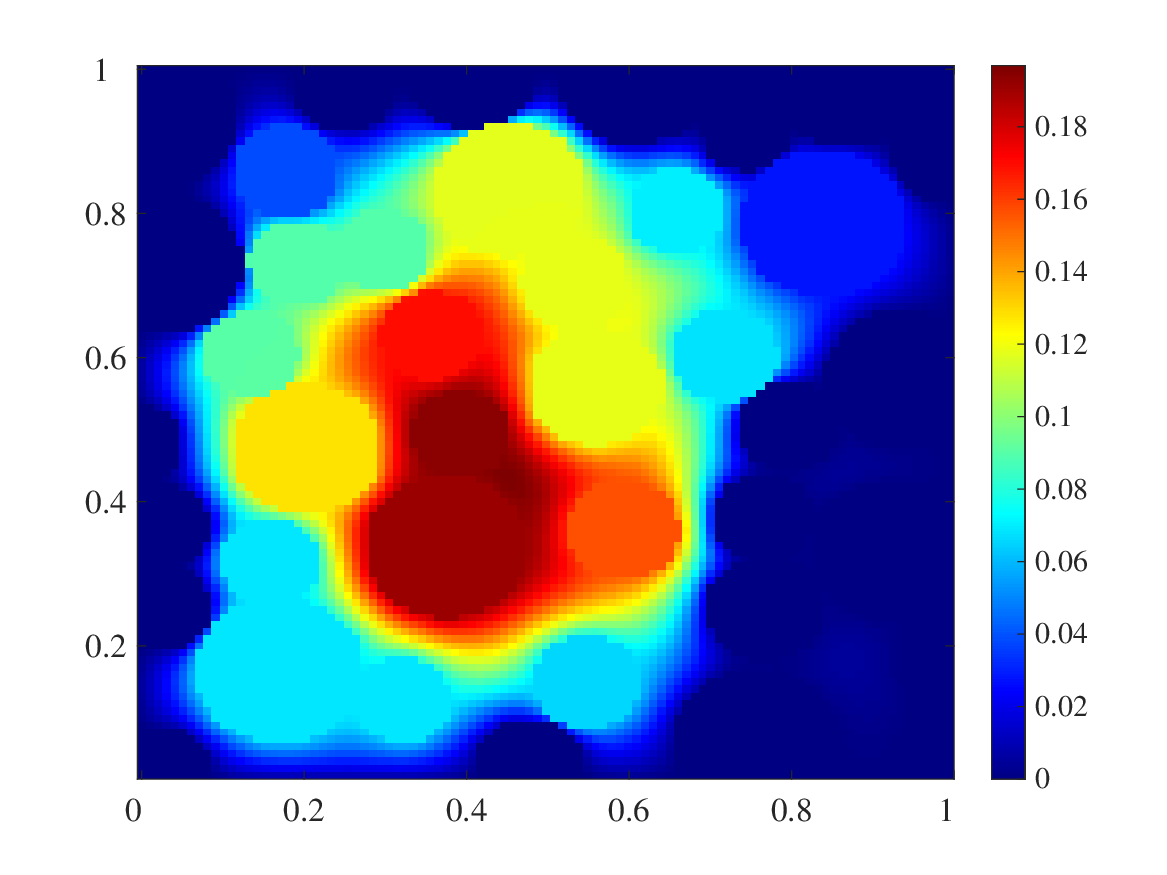}
  {\tiny(c)}\includegraphics[width=0.45\textwidth]{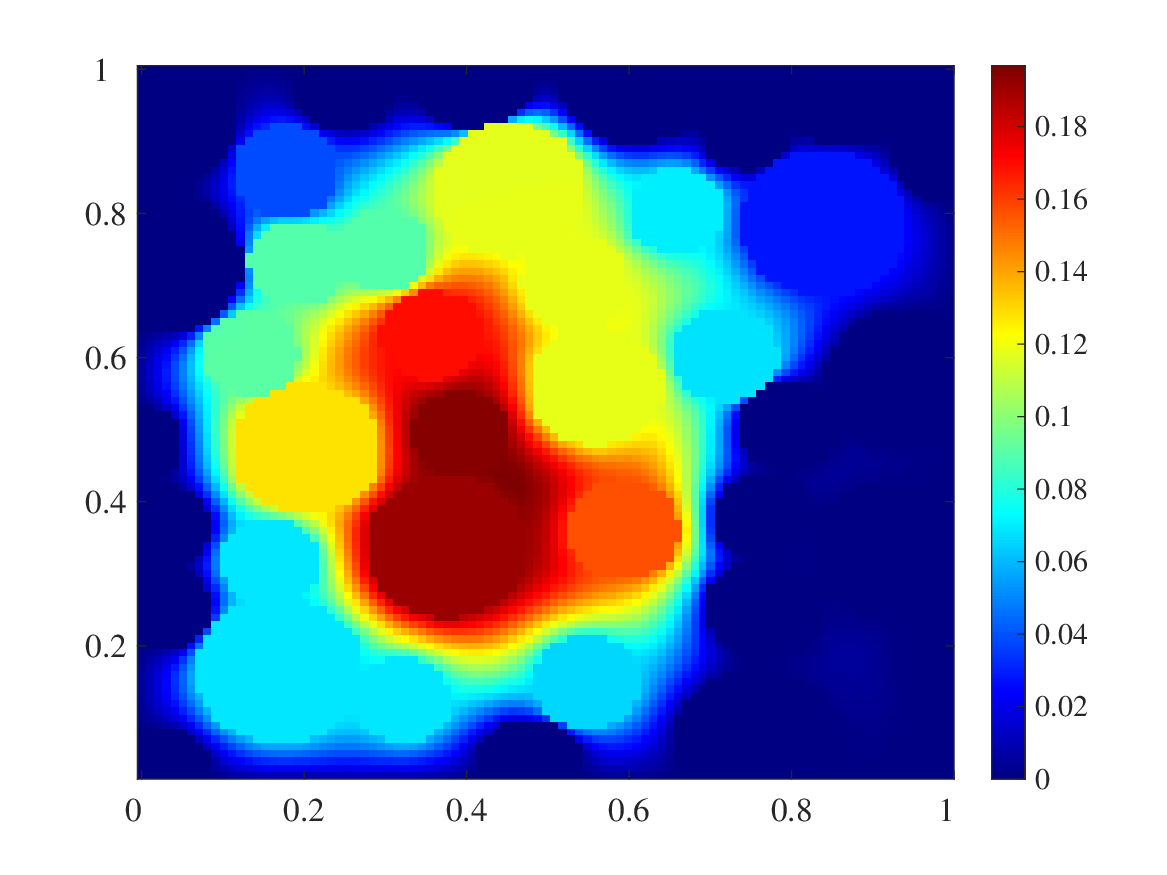}~
  {\tiny(f)}\includegraphics[width=0.45\textwidth]{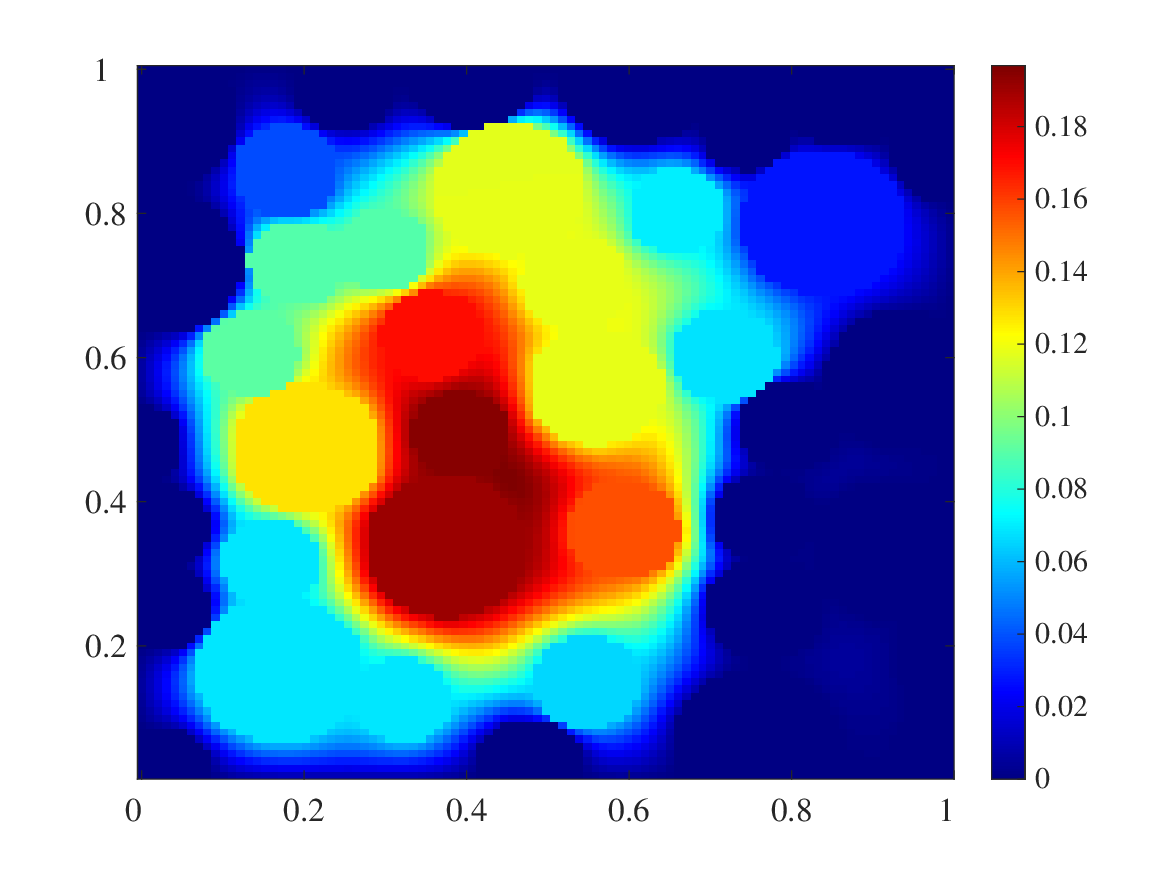}
  \caption{Contour plots of solutions: (a)the reference solution, (b)LOD,  (c) $LKSI\text{-}4$, (d)$LSSI\text{-}1$, (e)$LSSI\text{-}2$ and (f)$LSSI\text{-}4$.}
  \label{fig_ex1_solu}
\end{figure}

~\cref{fig_ex1_solu} displays solutions of several multiscale methods, including the LOD, LSSI and LKSI, where the number of oversampling layers is $m = 4$. Compared to the reference solution, all multiscale methods have captured the multiscale characteristics of the solution successfully and effectively. For further comparison, \cref{tab_ex1_compareMs} lists the energy error, $L^2$ error, degree of freedom (DoF), CPU time and number of local problems (NoLP) in multiscale methods. With an equivalent degree of freedom, both the LSSI and LKSI exhibit exceptional accuracy, and as the number of iteration steps $n$ increases, the error of the LSSI decays. \cref{tab_ex1_KSMsFEM} lists the results we are focusing on for the LKSI with different numbers of iteration steps (NoIS) $n$. As the number of iteration steps increases, there is a linear growth in the degrees of freedom, CPU time, and the number of local problems. Simultaneously, the energy error and $L^2$ error decrease at a decelerating rate.

\begin{table}[htbp]
  \centering
  \begin{tabular}{|l|c|c|c|c|c|}
  \hline
  Multiscale method                 & Energy  error & $L^2$ error   & DoF & CPU time (s) & NoLP \\ \hline
  $LOD $                              & 1.6780E-01    & 3.7891E-02 & 400 & 7.75     & 400                      \\ \hline
  $LSSI\text{-}1$  & 1.8494E-02      & 1.0610E-03 & 400 & 4.46     & 400                      \\ \hline
  $LSSI\text{-}2$  & 1.3449E-02      & 5.5521E-04 & 400 & 6.99     & 800                      \\ \hline
  $LSSI\text{-}4$  & 1.2695E-02      & 5.3744E-04 & 400 & 10.97    & 1600                     \\ \hline
  $LKSI\text{-}4$ & 1.1997E-02      & 5.3476E-04 & 400 & 9.88     & 400                      \\ \hline
  \end{tabular}
  \caption{Comparison of different multiscale methods in terms of the energy error, $L^2$ error, degree of freedom(DoF), CPU time and number of local problems(NoLP).}
  \label{tab_ex1_compareMs}
\end{table}

  \begin{table}[htbp]
    \centering
    \begin{tabular}{|c|c|c|c|c|c|}
    \hline
    NoIS & Energy error & $L^2$ error   & DoF & CPU time (s) & NoLP \\ \hline
    $LKSI\text{-}1$    & 2.9165E-02   & 3.3226E-03 & 100 & 2.53     & 100  \\ \hline
    $LKSI\text{-}2$    & 1.6337E-02   & 1.1029E-03 & 200 & 5.16     & 200  \\ \hline
    $LKSI\text{-}3$    & 1.2730E-02   & 6.6714E-04 & 300 & 7.35     & 300  \\ \hline
    $LKSI\text{-}4$    & 1.1997E-02   & 5.3476E-04 & 400 & 9.88     & 400  \\ \hline
    $LKSI\text{-}5$    & 1.2298E-02   & 5.1922E-04 & 500 & 11.67    & 500  \\ \hline
    \end{tabular}
    \caption{The energy error, $L^2$ error, degree of freedom(DoF), CPU time and number of local problems(NoLP) of the LKSI with different number of iteration steps (NoIS) $n$.}
    \label{tab_ex1_KSMsFEM}
  \end{table}

\begin{figure}[htbp]
    \centering
    {\tiny(a)}\includegraphics[width=0.45\textwidth]{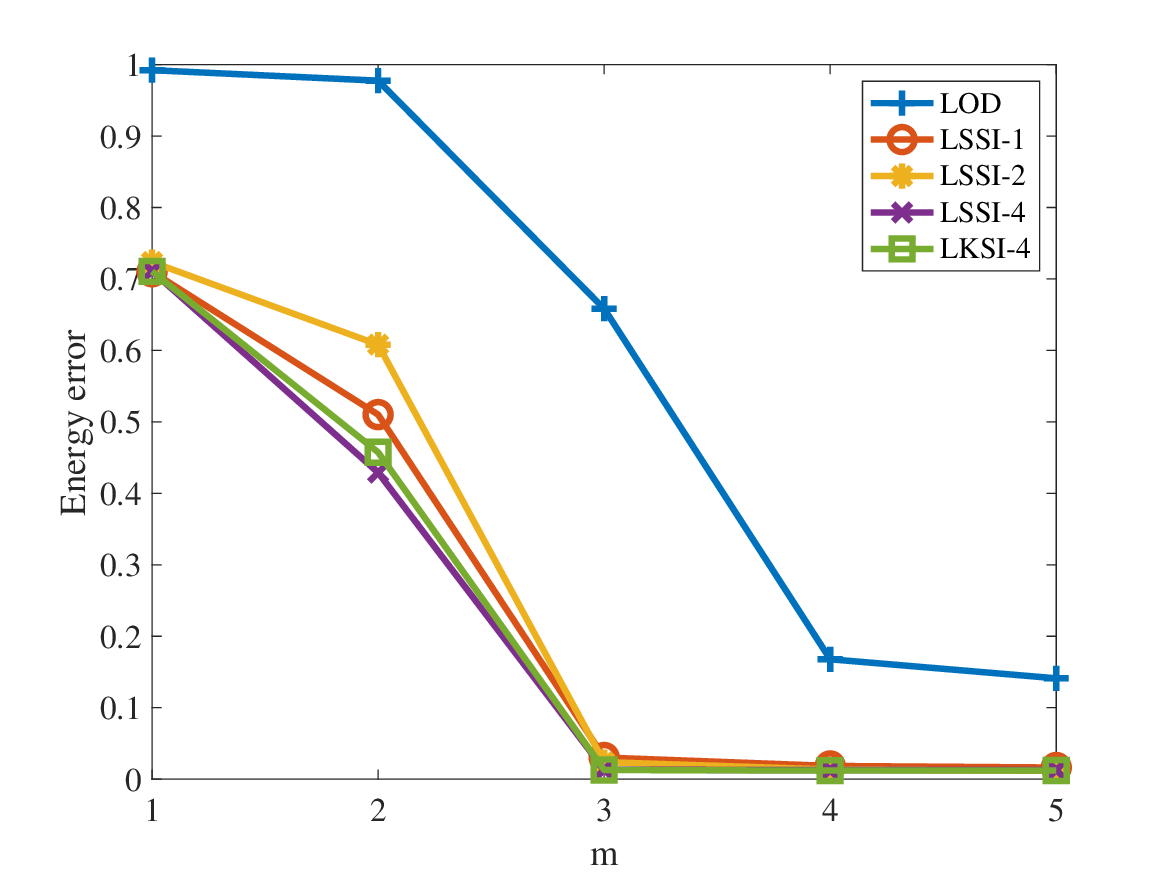}~
    {\tiny(b)}\includegraphics[width=0.45\textwidth]{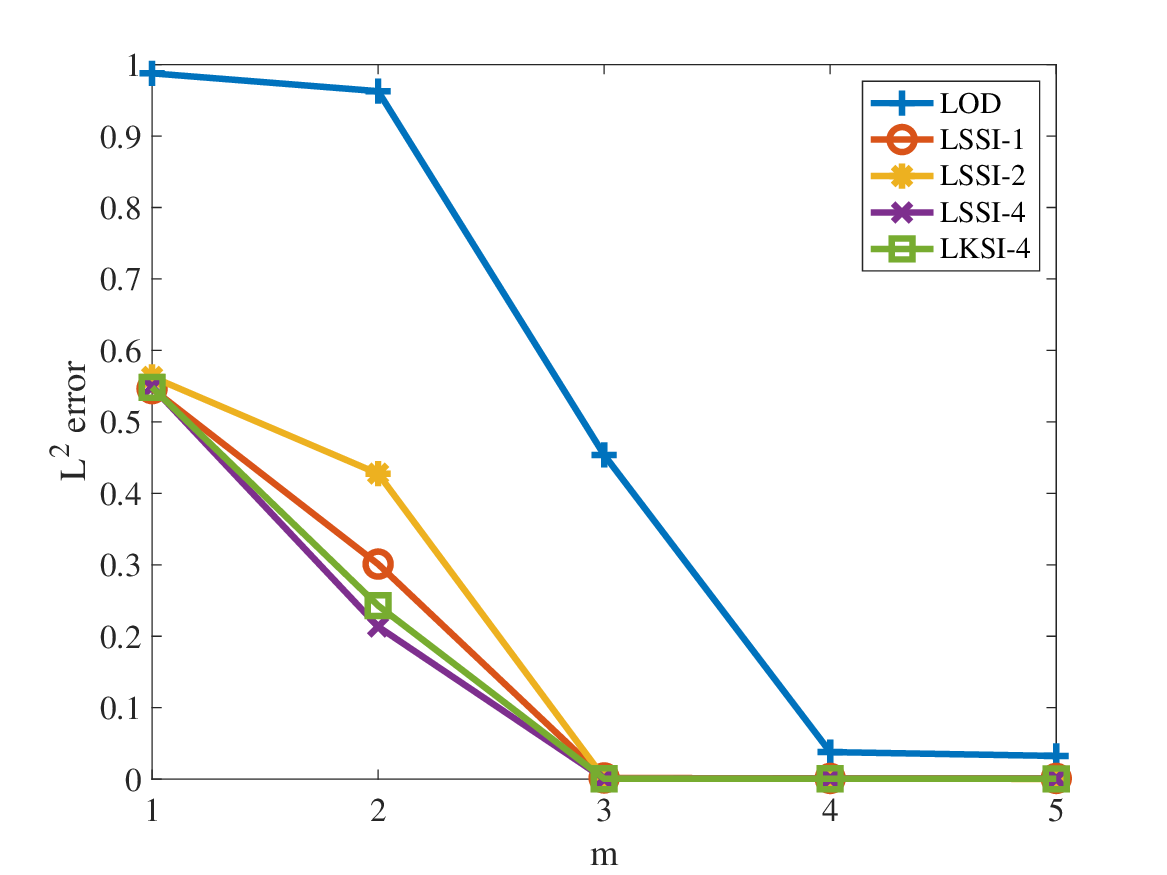}
    \caption{Relative errors of multiscale methods versus the number of oversampling layers $m$: (a) energy error and (b) $L^2$ error.}
    \label{fig_ex1_error_vs_m}
  \end{figure}

Among the many multiscale methods that employ oversampling techniques, the number of oversampling layers $m$ deserves significant consideration. When the coefficient $\kappa$ is high-contrast, the choice of $m$ in LOD is directly correlated with the contrast ${\kappa_{\max }}/{\kappa_{\min }}$. This is because the exponential decay rate of global basis functions is associated with the contrast \cite{RPD2021NHBSS}.
\cref{fig_ex1_error_vs_m} shows the relative errors of multiscale methods versus the number of oversampling layers $m$. In this numerical example, for the LSSI and LKSI, selecting $m = 3$ is sufficient, whereas for the LOD, $m$ needs to be at least 4.

\cref{tab_ex1_energy_contrast} and \cref{tab_ex1_L2_contrast} present the energy errors and $L^2$ errors of multiscale methods in relation to the contrast of $\kappa$, with a fixed number of oversampling layers $m=4$.
In the LOD, as the power exponent of contrast increases, there is a sharp rise in both the energy error and the $L^2$ error, reaching unacceptable levels. It is noteworthy that this issue can be mitigated by employing larger number of oversampling layers $m$. However, the relative errors are very consistent within the frameworks of LSSI and LKSI, even when there are big changes in the power exponent of contrast. Under appropriate conditions, we can argue that the relative errors of our proposed LSSI and LKSI are independent of contrast.
  \begin{table}[htbp]
    \begin{tabular}{|c|c|c|c|c|c|}
    \hline
    Contrast & $LOD$        & $LSSI\text{-}1$       & $LSSI\text{-}2$       & $LSSI\text{-}4$      & $LKSI\text{-}4$      \\ \hline
    1E+02        & 8.9027E-02 & 2.3644E-02 & 1.9436E-02 & 1.1528E-02 & 1.1330E-02 \\ \hline
    1E+03        & 1.2928E-01 & 2.0283E-02 & 1.5518E-02 & 1.2251E-02 & 1.1951E-02 \\ \hline
    1E+04        & 1.6780E-01 & 1.8494E-02 & 1.3449E-02 & 1.2695E-02 & 1.1997E-02 \\ \hline
    1E+05        & 3.0026E-01 & 1.7854E-02 & 1.2740E-02 & 1.2738E-02 & 1.1955E-02 \\ \hline
    1E+06        & 6.4146E-01 & 1.7532E-02 & 1.2615E-02 & 1.2621E-02 & 1.1906E-02 \\ \hline
    1E+07        & 9.0074E-01 & 1.7396E-02 & 1.2557E-02 & 1.5008E-02 & 1.1891E-02 \\ \hline
    \end{tabular}
    \caption{Energy errors of multiscale methods versus the contrast of $\kappa$, where the number of oversampling layers $m=4$.}
    \label{tab_ex1_energy_contrast}
    \end{table}

    \begin{table}[htbp]
      \begin{tabular}{|c|c|c|c|c|c|}
      \hline
      Contrast & $LOD$        & $LSSI\text{-}1$       & $LSSI\text{-}2$       & $LSSI\text{-}4$      & $LKSI\text{-}4$       \\ \hline
      1E+02        & 1.3417E-02 & 1.2870E-03 & 9.8800E-04 & 4.5969E-04 & 4.9329E-04 \\ \hline
      1E+03        & 2.6269E-02 & 1.1496E-03 & 6.7508E-04 & 4.9034E-04 & 5.3315E-04 \\ \hline
      1E+04        & 3.7892E-02 & 1.0610E-03 & 5.5521E-04 & 5.3744E-04 & 5.3476E-04 \\ \hline
      1E+05        & 9.3800E-02 & 1.0351E-03 & 5.1829E-04 & 5.5103E-04 & 5.2758E-04 \\ \hline
      1E+06        & 4.1972E-01 & 1.0314E-03 & 5.2722E-04 & 5.3654E-04 & 5.2040E-04 \\ \hline
      1E+07        & 8.2493E-01 & 1.0340E-03 & 5.2610E-04 & 6.3565E-04 & 5.1874E-04 \\ \hline
      \end{tabular}
      \caption{$L^2$ errors of multiscale methods versus the contrast of $\kappa$, where the number of oversampling layers $m=4$.}
      \label{tab_ex1_L2_contrast}
      \end{table}

\begin{figure}[htbp]
  \centering
  {\tiny(a)}\includegraphics[width=0.45\textwidth]{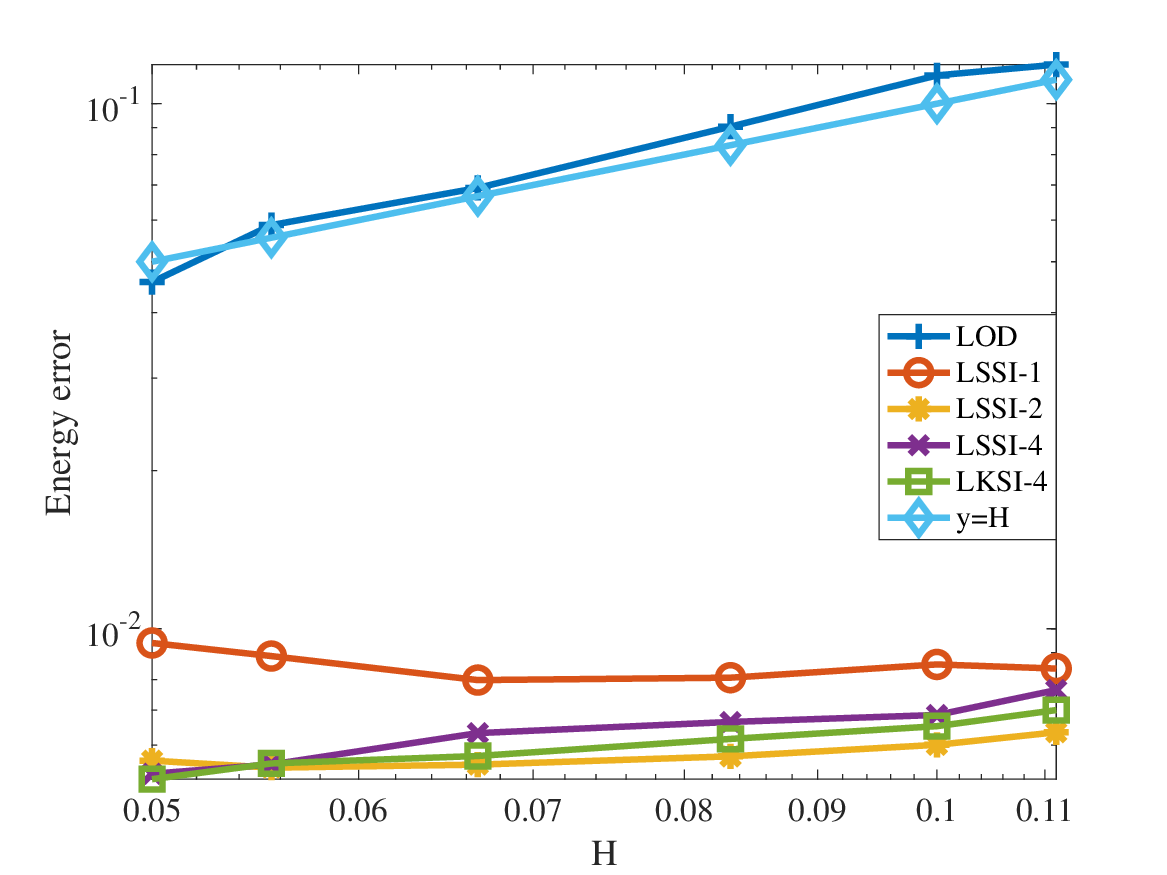}~
  {\tiny(b)}\includegraphics[width=0.45\textwidth]{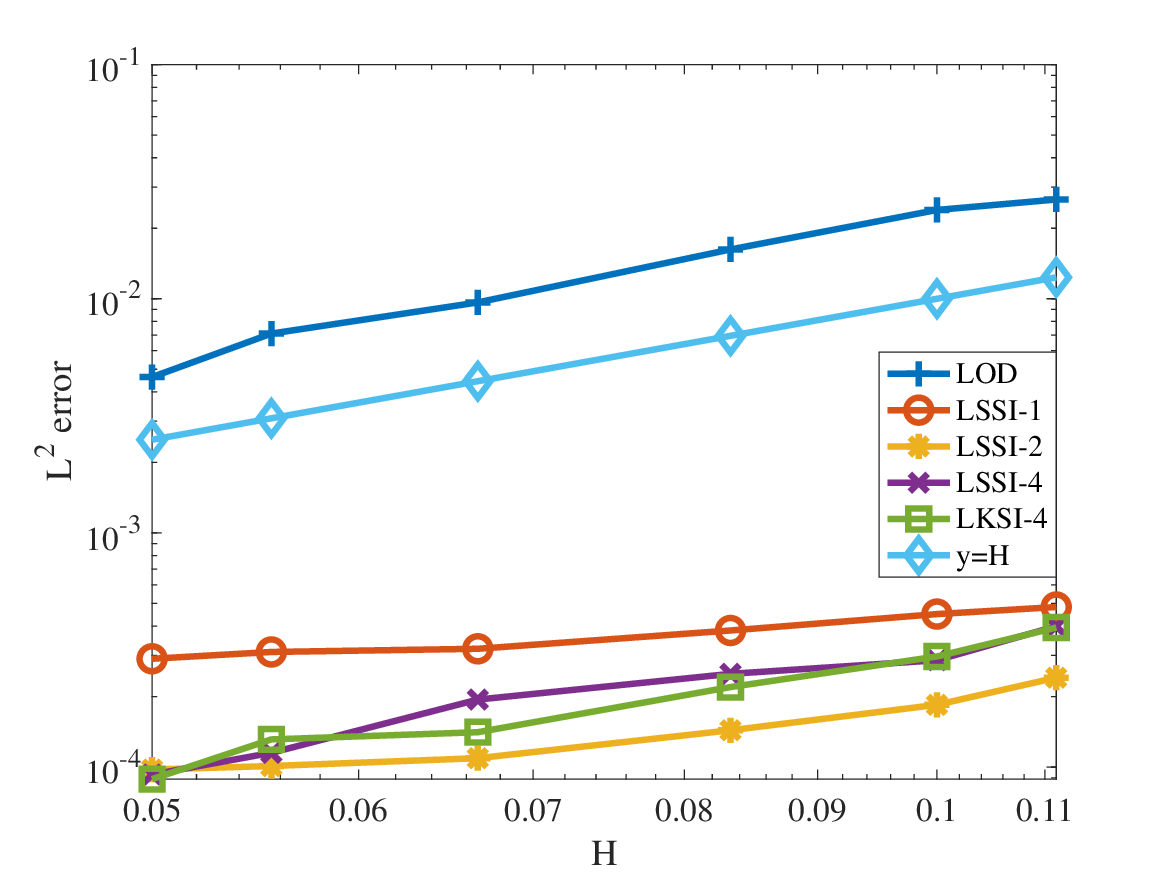}
  \caption{Relative errors of multiscale methods versus the coarse mesh size $H$: (a) energy error and (b) $L^2$ error.}
  \label{fig_ex1_error_vs_H}
\end{figure}

\cref{fig_ex1_error_vs_H} shows the relative errors of multiscale methods versus coarse mesh size $H$, where the fine mesh size is $h=1/180$ and the number of oversampling layers is $m = \lceil 2 \log (1 / H)\rceil$. For the LOD, the convergence rate of relative energy error is 1, and that of the relative $L^2$ error is 2, under an appropriate number of oversampling layers $m$, just as concluded in \cite{AD_LOD_2014}. Despite the fact that the convergence rates of relative errors in our proposed methods are lower than those of the LOD, the values themselves are significantly lower.

\subsection{Discussion of the channel length} \label{sec_num_ex2}
In multiscale problems with high-contrast coefficients, the focal challenge revolves around dealing with contrast, exemplified by the investigation of methods that remain independent of contrast variations. Interestingly, the geometric features of the large coefficient regions can also affect the results of multiscale methods. The influence of geometric features is very complex, and this article only focuses on channel length in simple cases.
\begin{figure}[htbp]
  \centering
  {\tiny(a)}\includegraphics[width=0.3\textwidth]{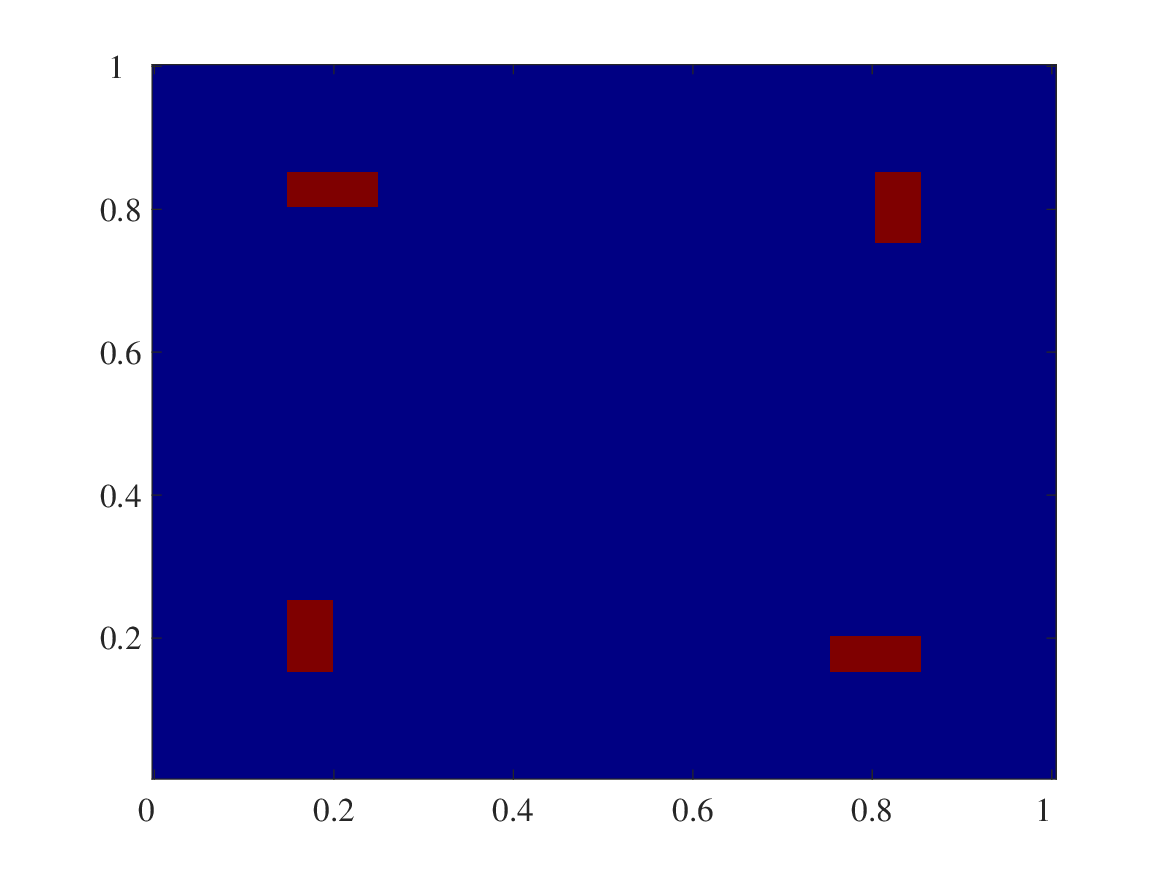}~
  {\tiny(b)}\includegraphics[width=0.3\textwidth]{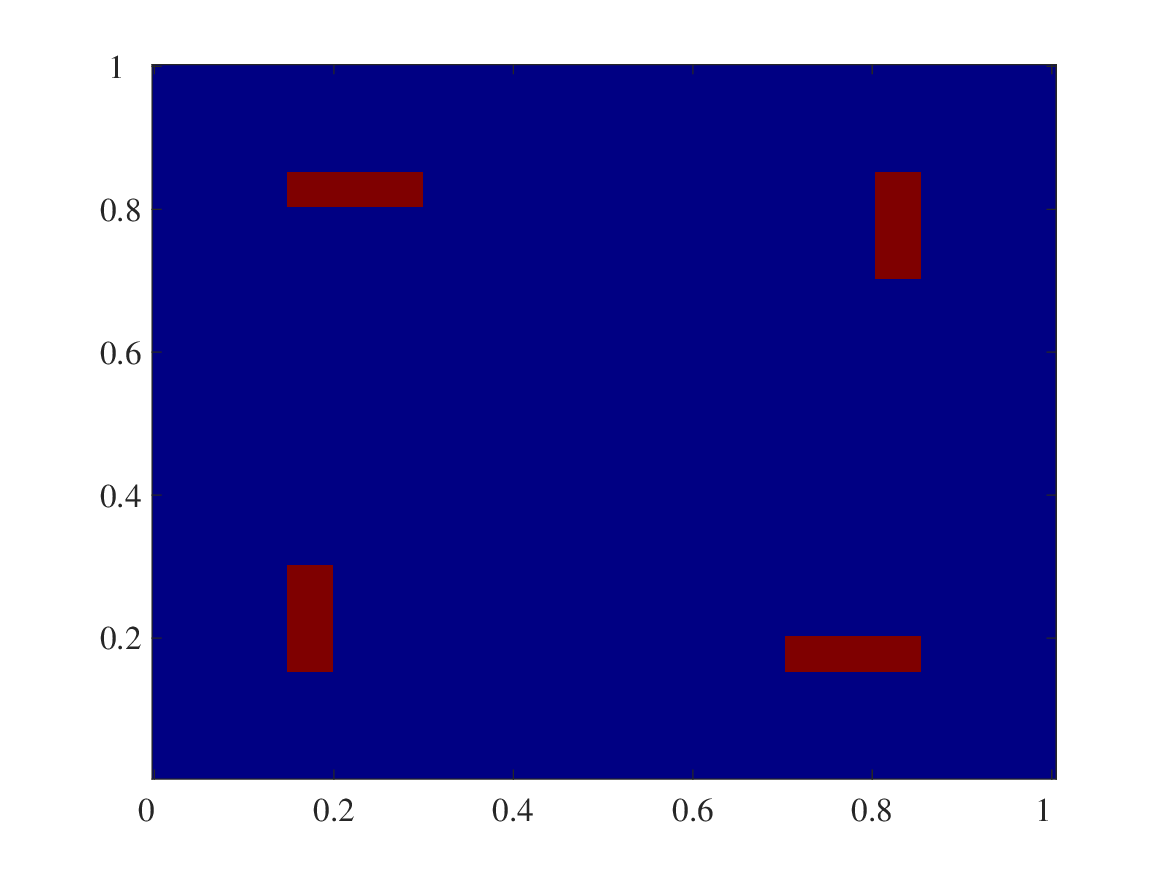}~
  {\tiny(c)}\includegraphics[width=0.3\textwidth]{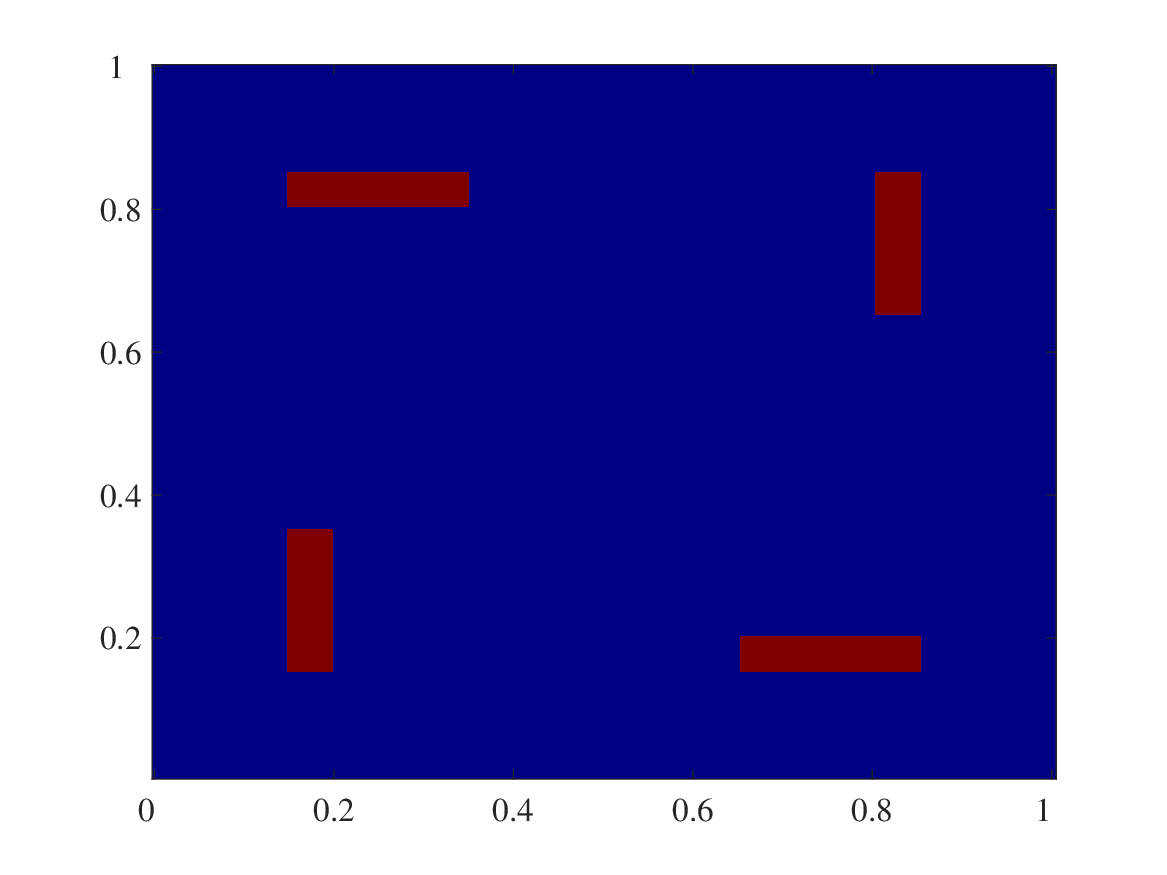}
  {\tiny(d)}\includegraphics[width=0.3\textwidth]{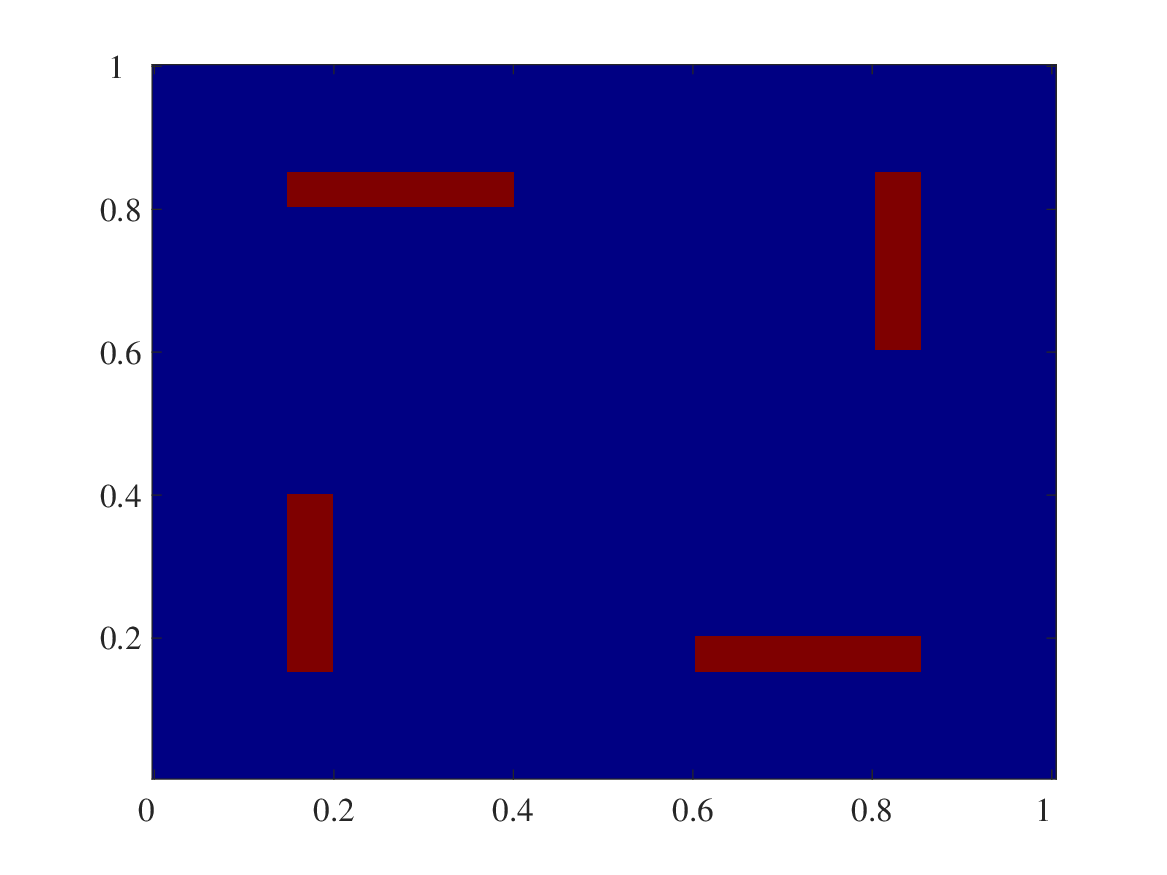}~
  {\tiny(e)}\includegraphics[width=0.3\textwidth]{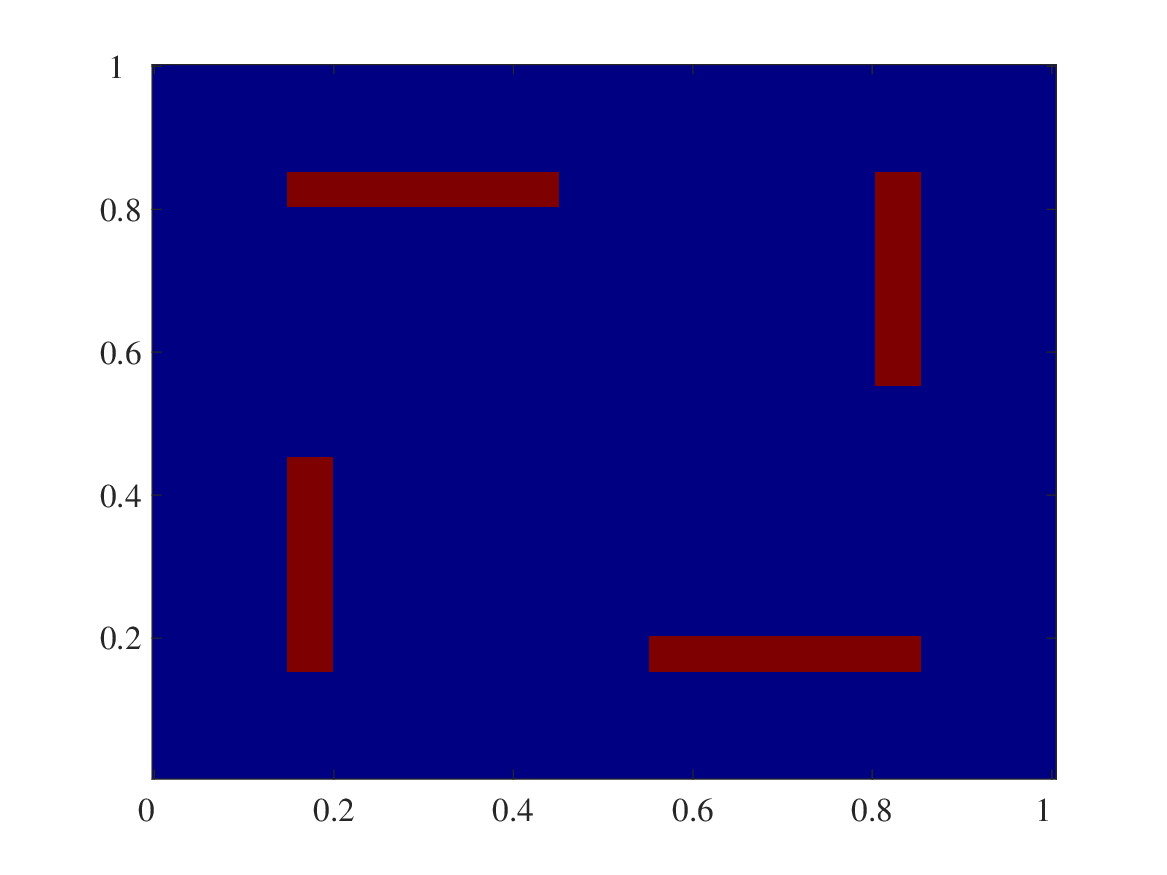}~
  {\tiny(f)}\includegraphics[width=0.3\textwidth]{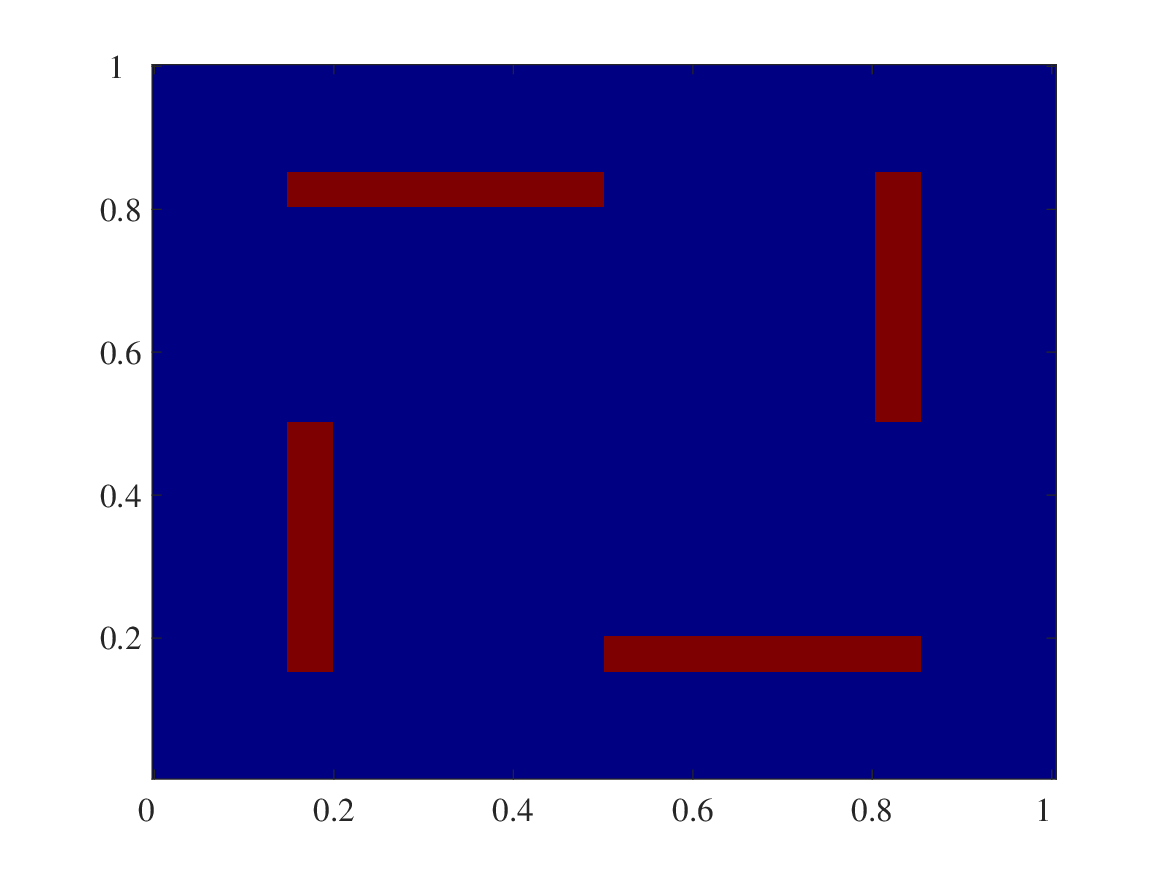}
  {\tiny(g)}\includegraphics[width=0.3\textwidth]{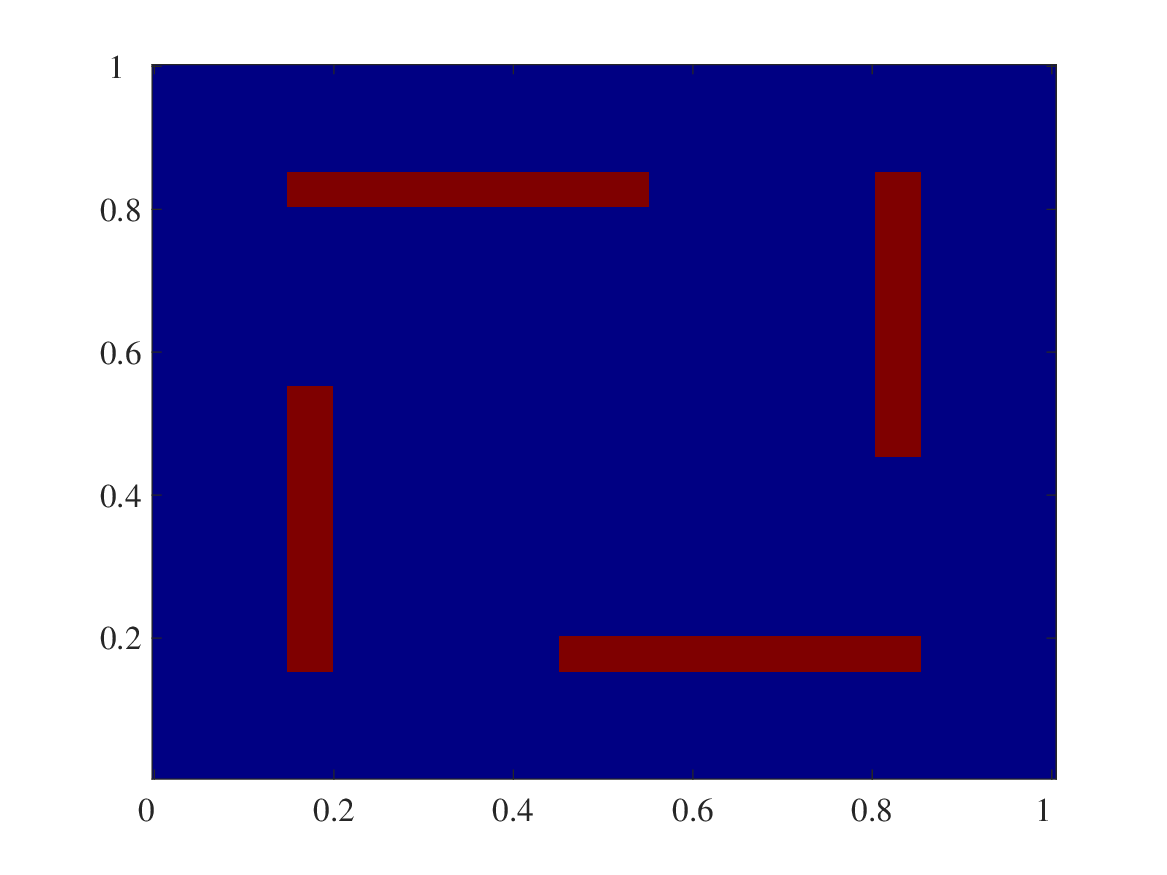}~
  {\tiny(h)}\includegraphics[width=0.3\textwidth]{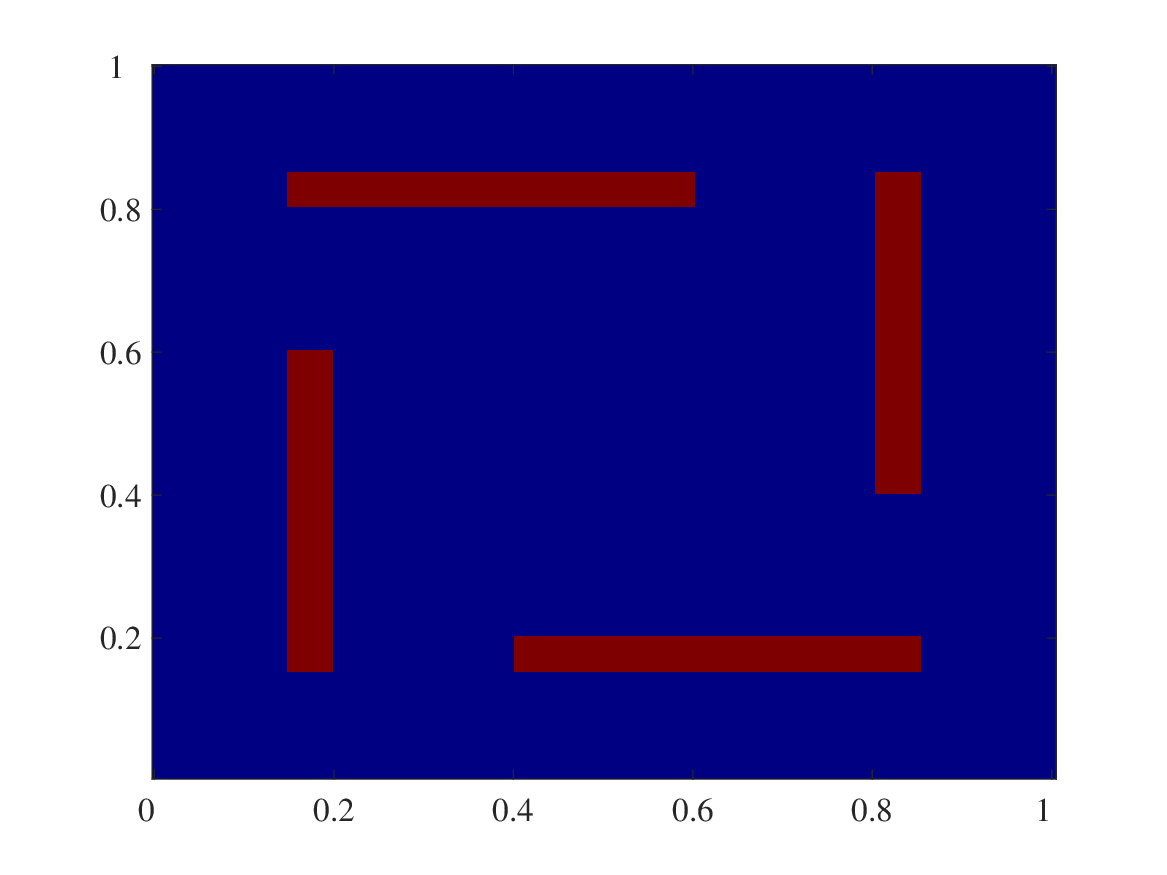}~
  {\tiny(i)}\includegraphics[width=0.3\textwidth]{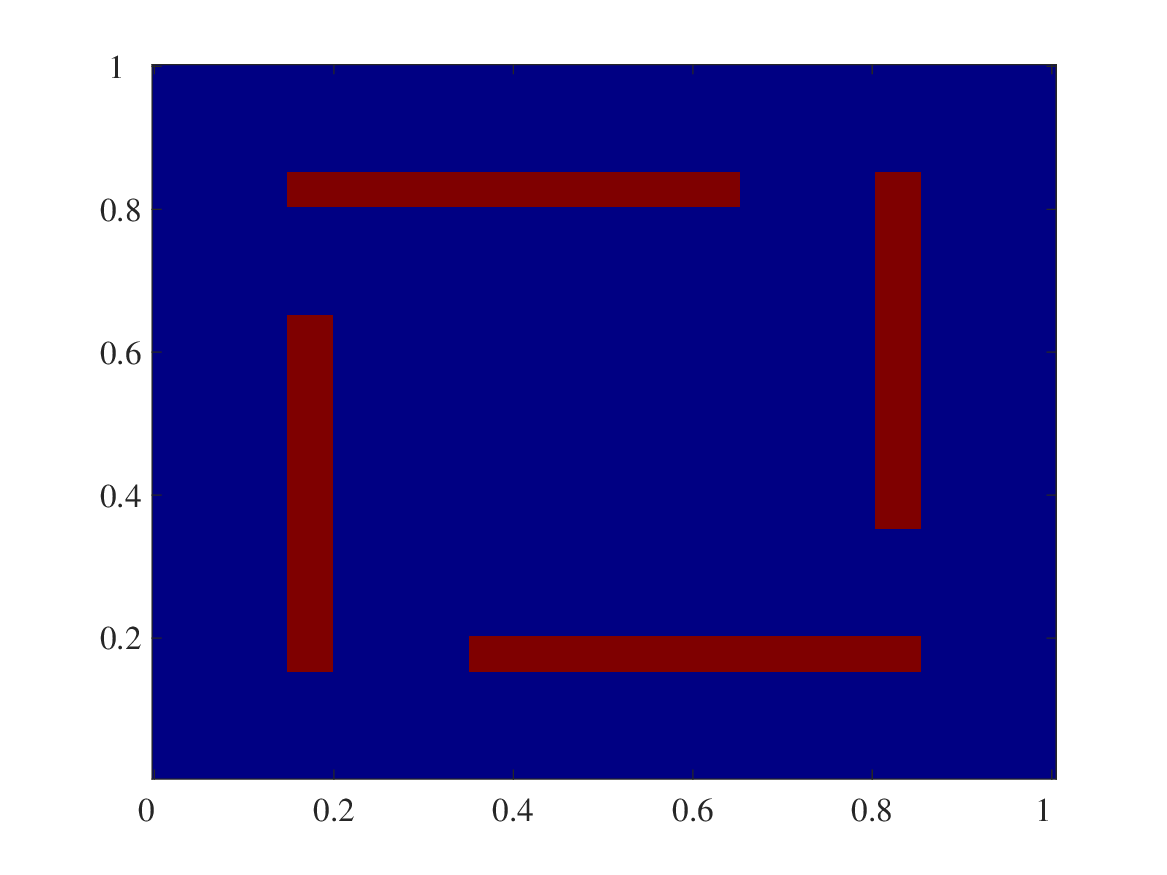}
  \caption{Multiscale coefficients with different channel lengths.}
  \label{fig_ex2_kappa}
\end{figure}

\begin{figure}[htbp]
  \centering
  {\tiny(a)}\includegraphics[width=0.45\textwidth]{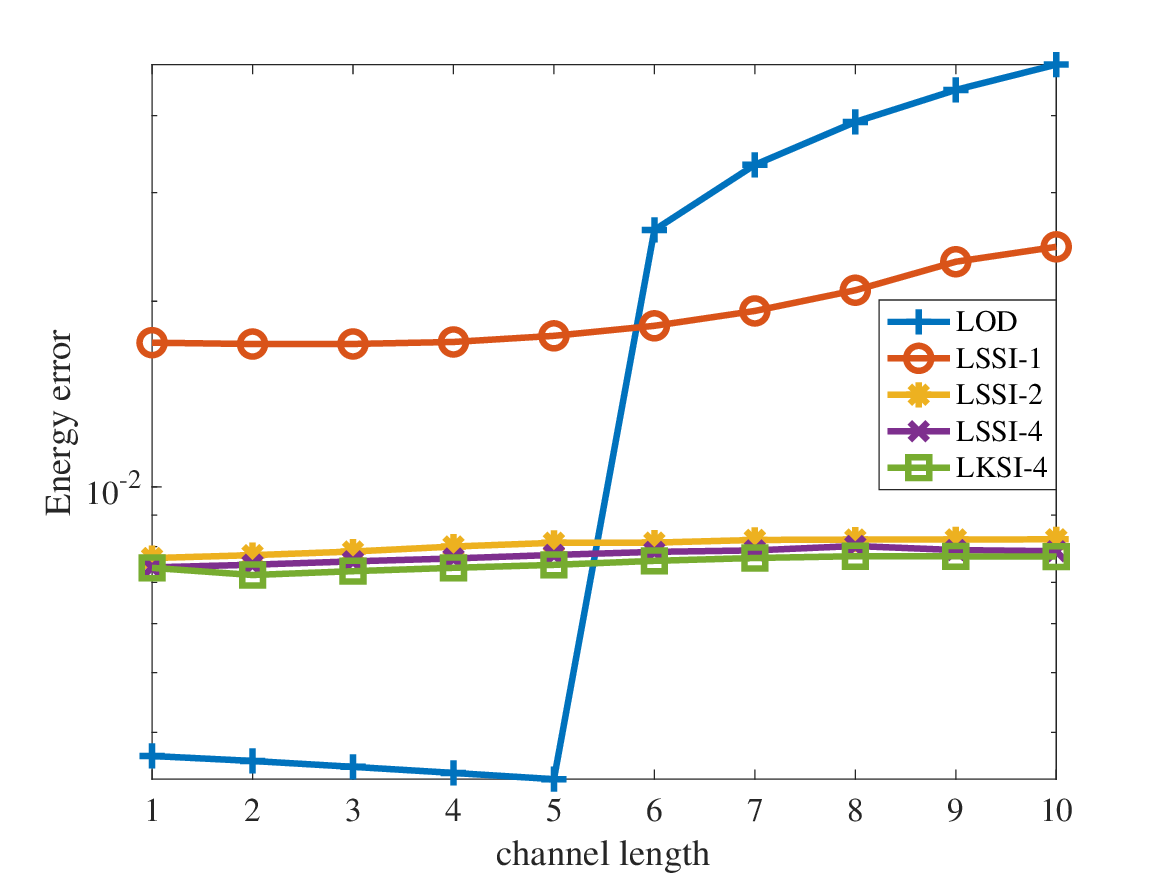}~
  {\tiny(b)}\includegraphics[width=0.45\textwidth]{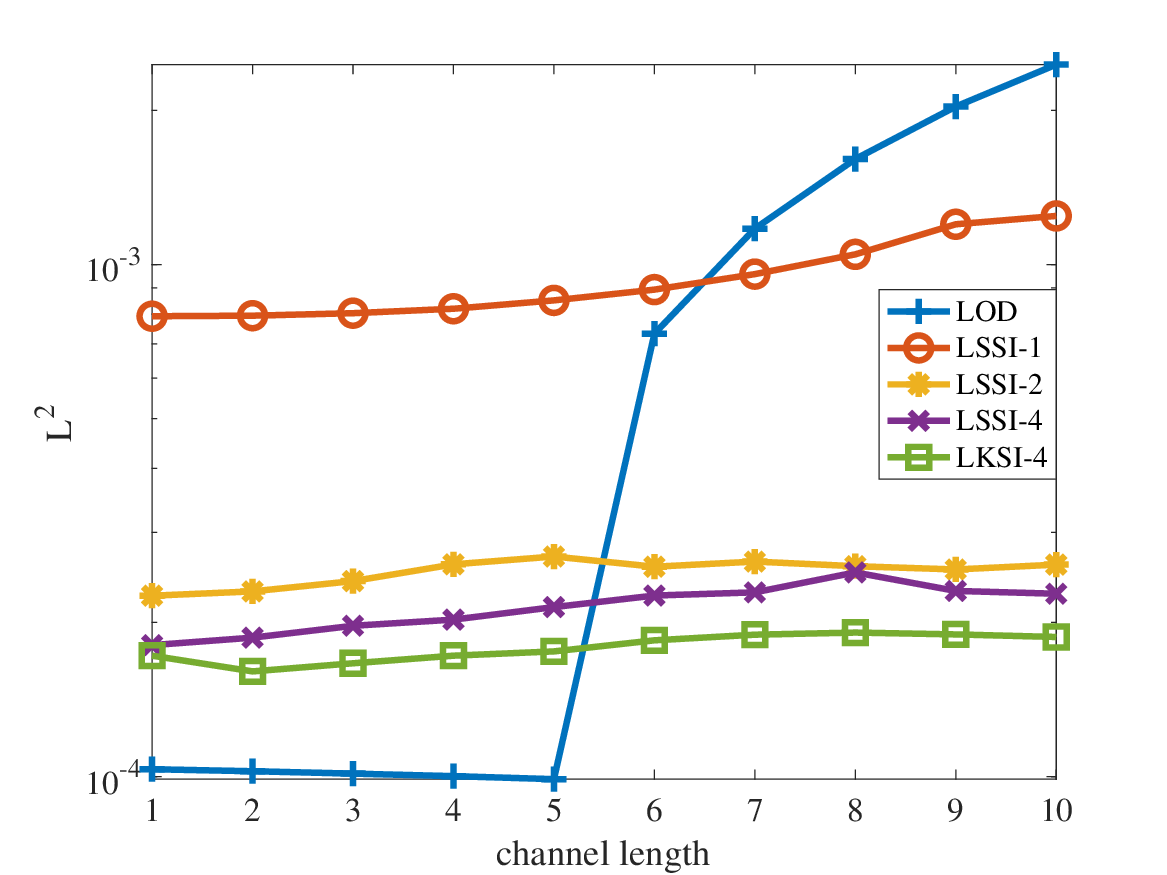}
  \caption{Relative errors of multiscale methods versus the channel length: (a) energy error and (b) $L^2$ error.}
  \label{fig_ex2_error_vs_cl}
\end{figure}

We also consider problem \cref{eq_diffusion}, and \cref{fig_ex2_kappa} shows multiscale coefficients $\kappa$ with different channel lengths (from $2H$ to $10H$). The fine mesh size is $h=1/200$, the coarse mesh size is $H=1/20$, the number of oversampling layers is $m=5$ and the contrast is $10^4$. Other problem settings are the same as in the previous numerical example. \cref{fig_ex2_error_vs_cl} shows the relative errors of multiscale methods versus the channel length. When the channel length is no higher than 5, the relative errors of the LOD are significantly better than those of our proposed methods. But when the channel length is greater than 5, the errors of the LOD increase sharply, while those of our proposed method increase slightly. This shows that our proposed method has stronger stability for long-channel cases.

\subsection{Elasticity problem}
The methods we introduced are not limited to diffusion problems; in fact, they hold for general positive definite operators. We consider an elasticity problem
\begin{equation} \label{elasticity}
      -\nabla\cdot\mathbf{\sigma}\big( \mathbf{u} \big)=\mathbf{f}, \text{ in } {\Omega}, \\
\end{equation}
where $\Omega = [0,1]^2$ and $\mathbf{f} =[ \sin(\pi x) \sin(\pi y), 1 ] $. The stress-strain relationship is given by
\begin{equation}
    \mathbf{\sigma}\big( \mathbf{u} \big)=2 {\mu } \epsilon(u)+{\lambda }\nabla \cdot  \mathbf{u} I,
\end{equation}
where $\lambda >0$ and $\mu >0$ are the Lam$\acute{\text{e}}$ constants. The strain tensor $ \epsilon(\mathbf{u}) = (\epsilon_{ij}(\mathbf{u}))_{1 \leq i,j \leq 2}$ is defined by
\begin{equation}
    \epsilon( \mathbf{u} )=\frac{1}{2}\left(\nabla \mathbf{u} +\nabla \mathbf{u}^{T}\right).
\end{equation}
The Lam$\acute{\text{e}}$ constants $\lambda$ and $\mu$ are the same as $\kappa$ in \cref{sec_NumExam_ex1}. The fine mesh size is $h=1/100$, and the coarse mesh size is $H=1/10$. In the LOD and LSSI, we select $8$ bilinear functions on $K_i$ as the initial basis functions $\left\{\phi_i^{j,0} \right\}_{j=1}^{8}$ for each subdomain $K_i^{m}$. In the LKSI, we select piecewise constant on $K_i$ as the initial basis function $\psi_i^{0} $ for each subdomain $K_i^{m}$.
\begin{figure}[htbp]
  \centering
  {\tiny(a)}\includegraphics[width=0.45\textwidth,height = 2.5cm]{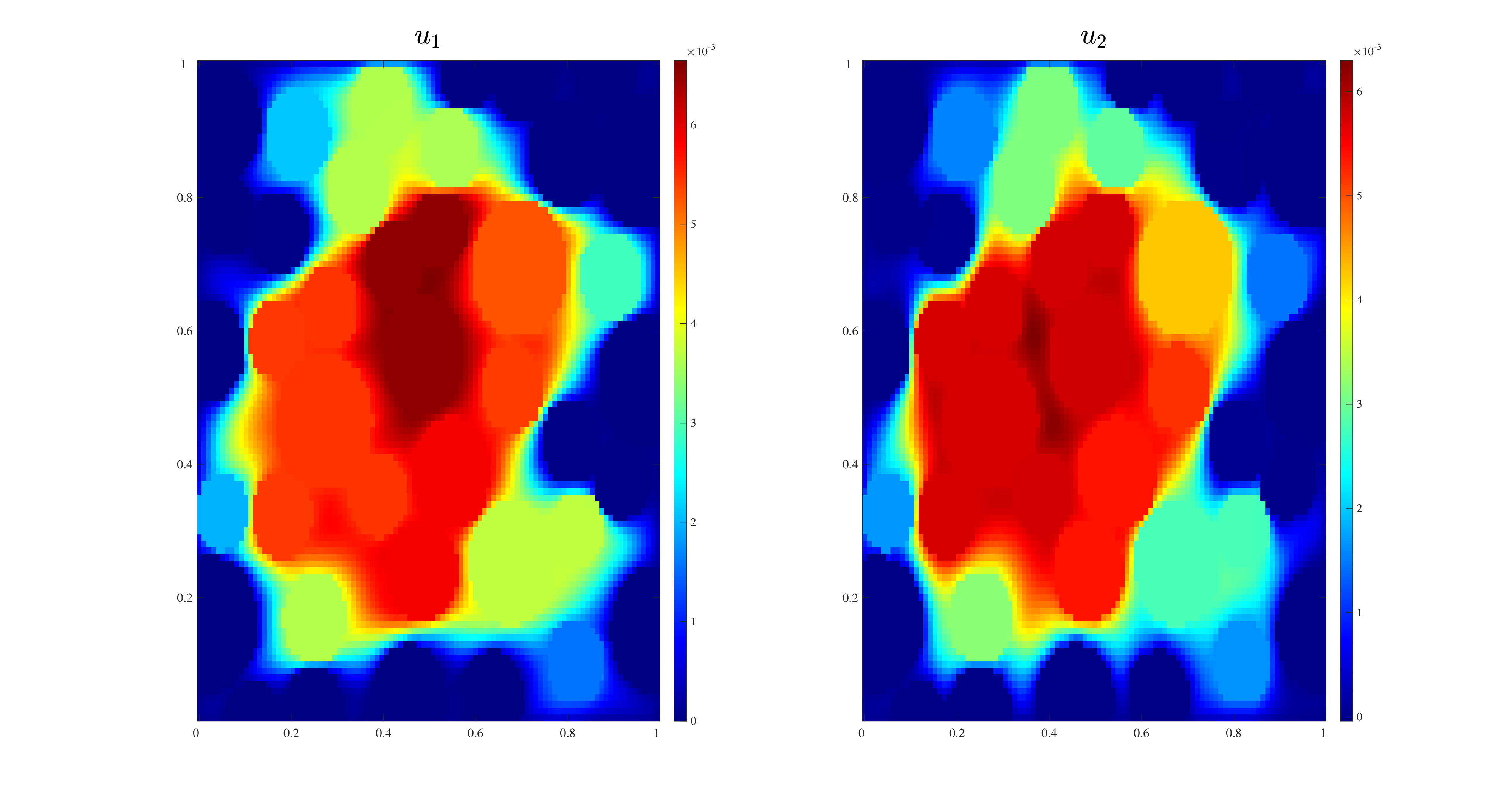}~
  {\tiny(d)}\includegraphics[width=0.45\textwidth,height = 2.5cm]{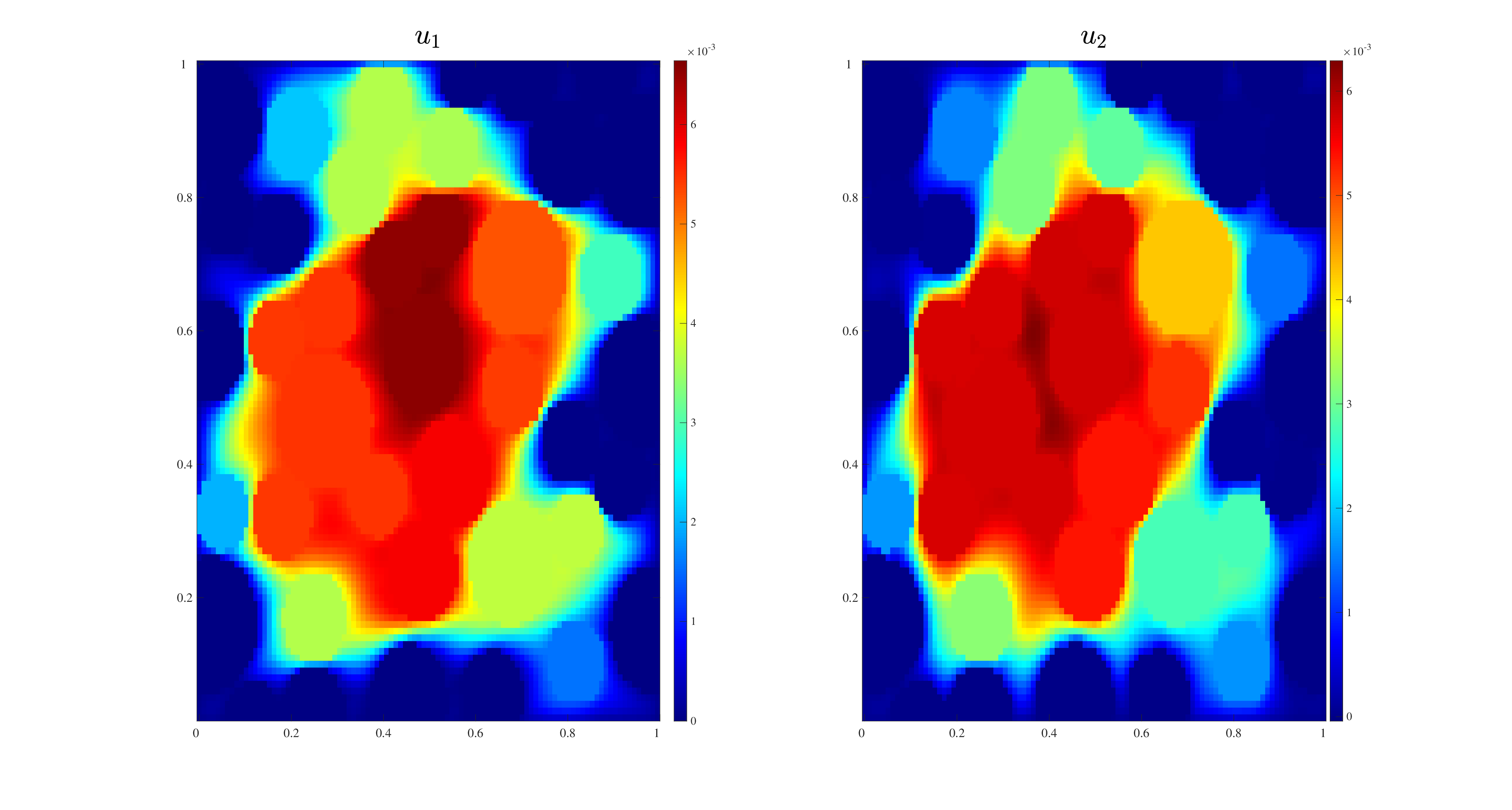}
  {\tiny(b)}\includegraphics[width=0.45\textwidth,height = 2.5cm]{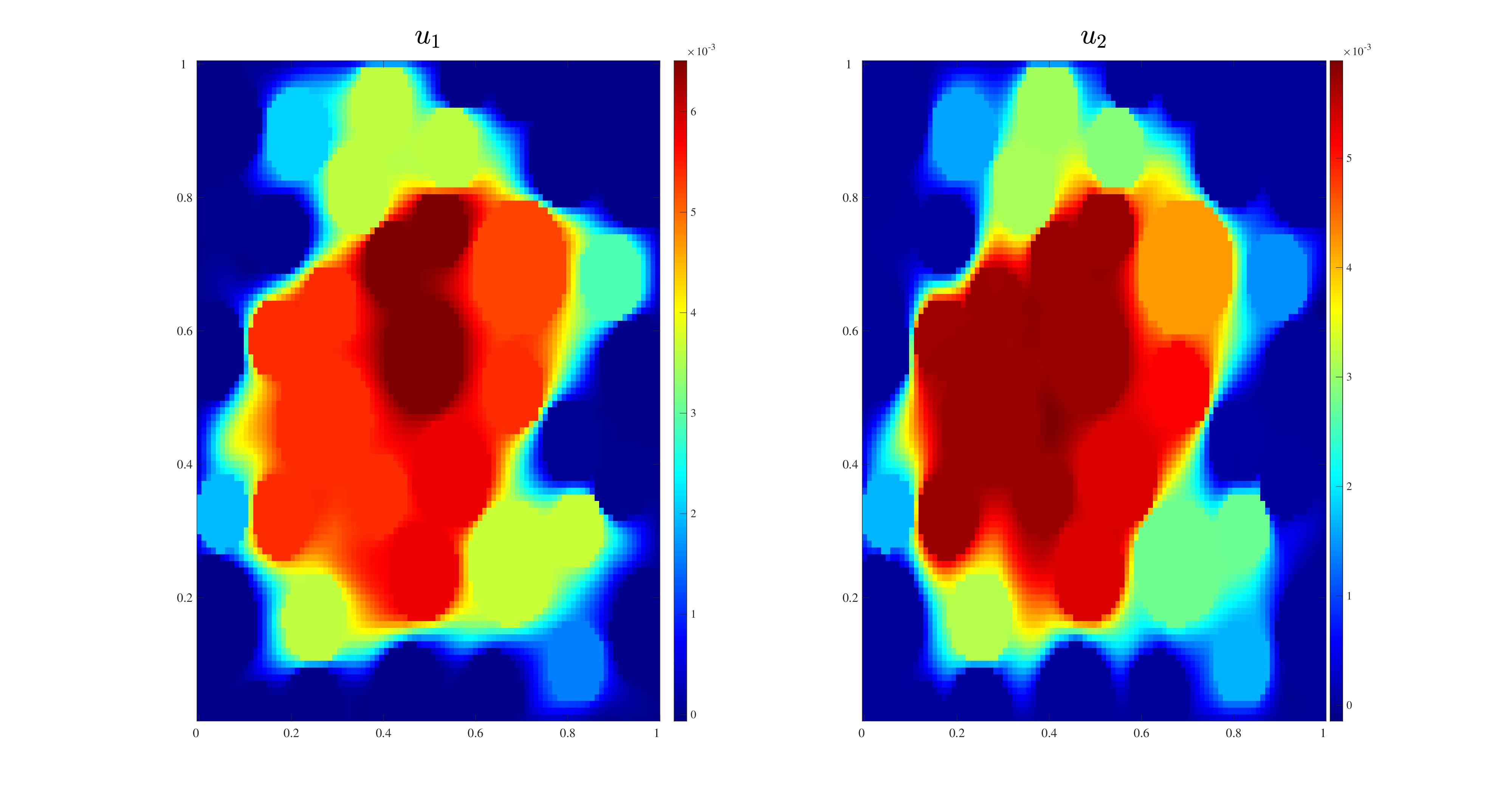}~
  {\tiny(e)}\includegraphics[width=0.45\textwidth,height = 2.5cm]{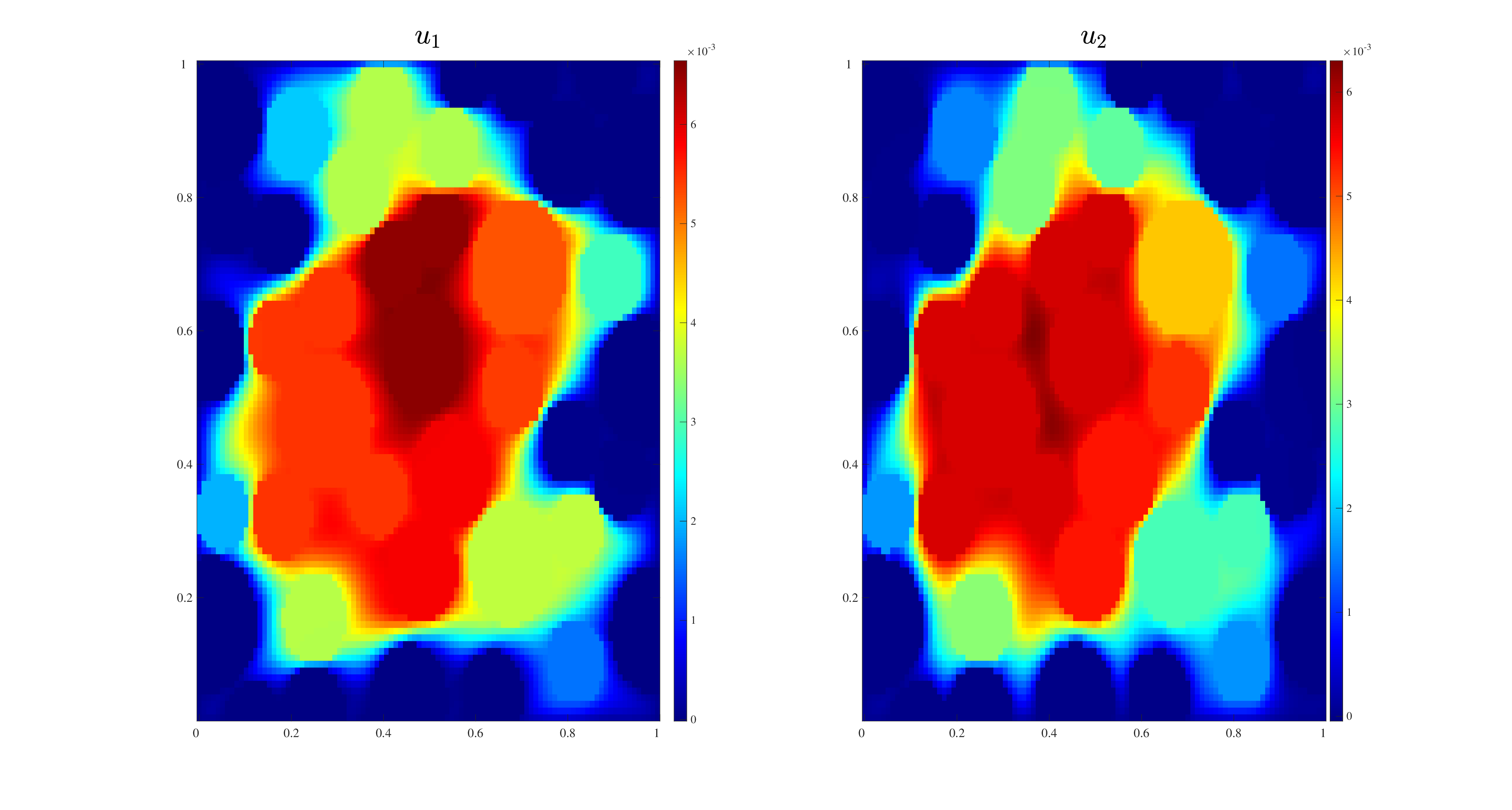}
  {\tiny(c)}\includegraphics[width=0.45\textwidth,height = 2.5cm]{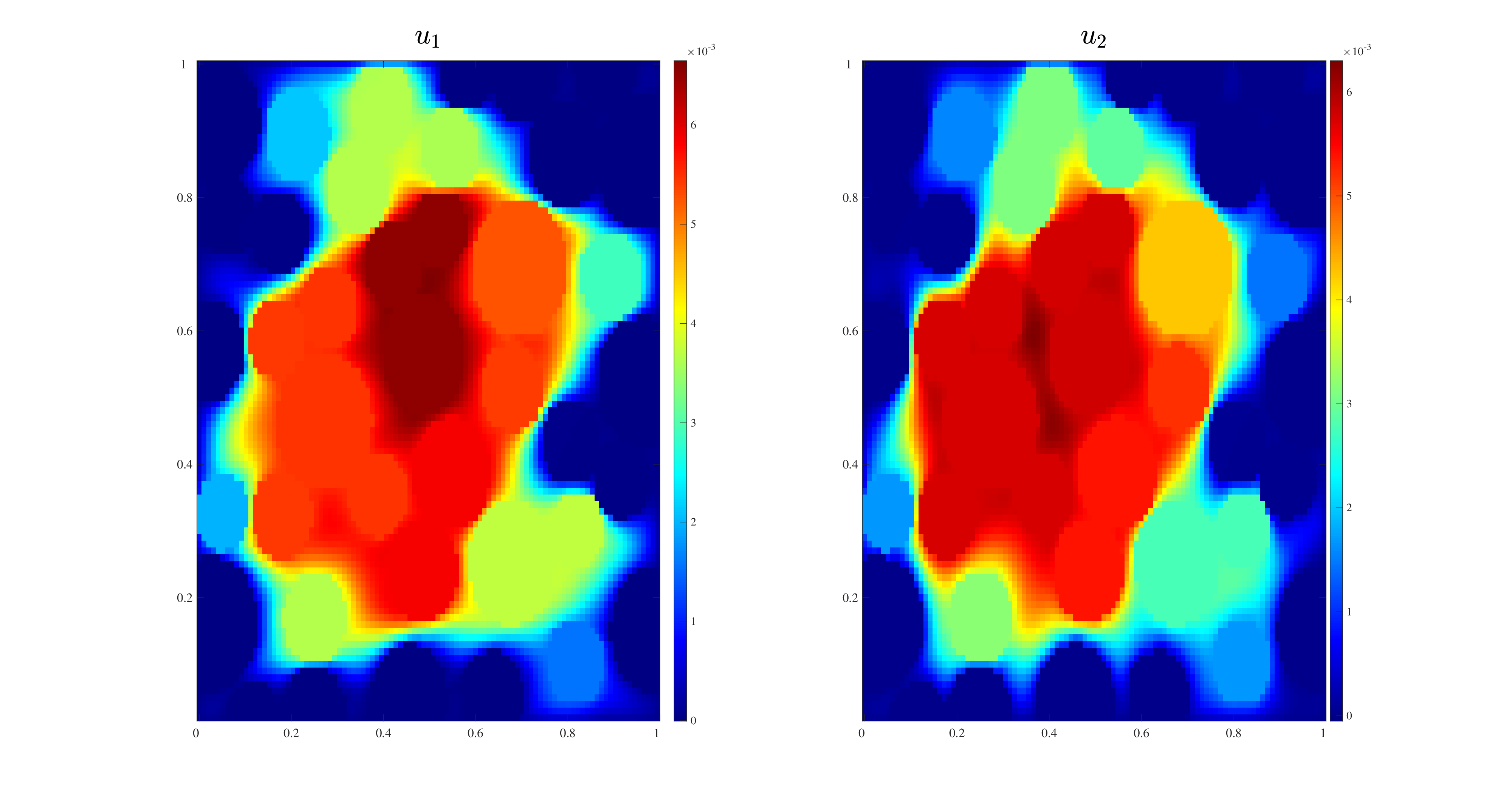}~
  {\tiny(f)}\includegraphics[width=0.45\textwidth,height = 2.5cm]{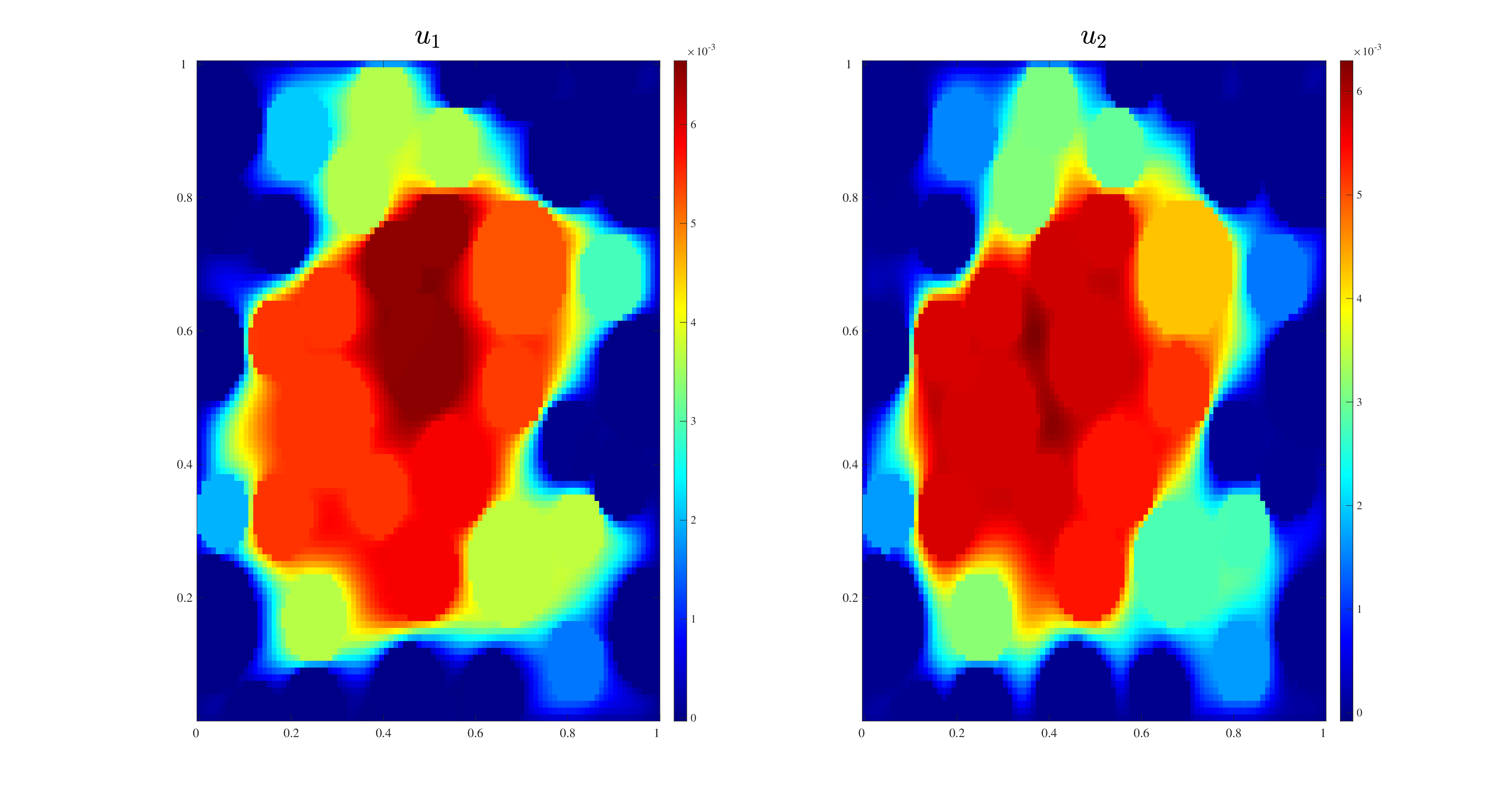}
  \caption{Contour plots of solutions: (a)the reference solution, (b)LOD,  (c) $LKSI\text{-}6$, (d)$LSSI\text{-}1$, (e)$LSSI\text{-}2$ and (f)$LSSI\text{-}8$. In each subfigure, the left is the first component of the deformation $\mathbf{u}$, and the right side is the second component.}
  \label{fig_ex3_solu}
\end{figure}

\begin{table}[htbp]
  \centering
  \begin{tabular}{|l|c|c|c|c|c|}
  \hline
  Multiscale method                 & Energy  error & $L^2$ error   & DoF & CPU time (s) & NoLP \\ \hline
  $LOD $                              & 1.8414E-01   & 4.5567E-02 & 800 & 49.62     & 800                      \\ \hline
  $LSSI\text{-}1$  & 1.7168E-02      & 1.3079E-03 & 800 & 18.30    & 800                      \\ \hline
  $LSSI\text{-}2$  & 1.0592E-02      & 8.1541E-04 & 800 & 27.10     & 1600                      \\ \hline
  $LSSI\text{-}8$  & 1.7368E-02      & 1.5965E-03 & 800 & 78.94    & 6400                     \\ \hline
  $LKSI\text{-}6$ & 9.7085E-03      & 9.7971E-04 & 600 & 42.54     & 600                      \\ \hline
  \end{tabular}
  \caption{Comparison of different multiscale methods in terms of the energy error, $L^2$ error, degree of freedom(DoF), CPU time and number of local problems(NoLP).}
  \label{tab_ex3_compareMs}
\end{table}

\begin{table}[htbp]
  \centering
  \begin{tabular}{|c|c|c|c|c|c|}
  \hline
  NoIS & Energy  error   & $L^2$ error   & DoF & CPU time & NoLP \\ \hline
  $LSSI\text{-}1$              & 1.7168E-02 & 1.3079E-03 & 800 & 18.30    & 800                      \\ \hline
  $LSSI\text{-}2$              & 1.0592E-02 & 8.1541E-04 & 800 & 27.10    & 1600                     \\ \hline
  $LSSI\text{-}3$              & 1.1081E-02 & 7.9146E-04 & 800 & 34.63    & 2400                     \\ \hline
  $LSSI\text{-}4$              & 1.2494E-02 & 9.3155E-04 & 800 & 43.46    & 3200                     \\ \hline
  $LSSI\text{-}5$              & 1.3972E-02 & 1.0856E-03 & 800 & 51.97    & 4000                     \\ \hline
  $LSSI\text{-}6$              & 1.5448E-02 & 1.2735E-03 & 800 & 61.36    & 4800                     \\ \hline
  $LSSI\text{-}7$              & 1.6625E-02 & 1.4493E-03 & 800 & 70.01    & 5600                     \\ \hline
  $LSSI\text{-}8$              & 1.7308E-02 & 1.5965E-03 & 800 & 78.94    & 6400                     \\ \hline
  \end{tabular}
  \caption{The energy error, $L^2$ error, degree of freedom(DoF), CPU time and number of local problems(NoLP) of the LSSI with different number of iteration steps (NoIS) $n$.}
  \label{tab_ex3_SMsFEM}
  \end{table}

  \begin{table}[htbp]
    \centering
    \begin{tabular}{|c|c|c|c|c|c|}
    \hline
    NoIS & Energy error   & $L^2$ error   & DoF & CPU time & NoLP \\ \hline
    $LKSI\text{-}1$              & 6.0448E-02 & 2.4621E-02 & 100 & 7.59     & 100                      \\ \hline
    $LKSI\text{-}2$             & 2.7115E-02 & 3.9845E-03 & 200 & 13.96    & 200                      \\ \hline
    $LKSI\text{-}3$            & 1.5114E-02 & 1.7718E-03 & 300 & 20.44    & 300                      \\ \hline
    $LKSI\text{-}4$             & 1.2045E-02 & 1.2045E-02 & 400 & 26.7     & 400                      \\ \hline
    $LKSI\text{-}5$             & 1.0695E-02 & 1.1438E-03 & 500 & 33.47    & 500                      \\ \hline
    $LKSI\text{-}6$            & 9.7085E-03 & 9.7971E-04 & 600 & 42.54    & 600                      \\ \hline
    \end{tabular}
    \caption{The energy error, $L^2$ error, degree of freedom(DoF), CPU time and number of local problems(NoLP) of the LKSI with different number of iteration steps (NoIS) $n$.}
    \label{tab_ex3_KSMsFEM}
    \end{table}

\cref{fig_ex3_solu} displays solutions of several multiscale methods, where the number of oversampling layers is $m = 4$. \cref{tab_ex3_compareMs} lists the energy error, $L^2$ error, degree of freedom (DoF), CPU time and number of local problems (NoLP) in various multiscale methods. The results obtained are almost consistent with those obtained in the previous numerical examples. \cref{tab_ex3_SMsFEM} and \cref{tab_ex3_KSMsFEM} list the results of the LSSI and LKSI that we focus on with different number of iteration steps (NoIS) $n$. For the LSSI, the increase in the number of iteration steps means that the obtained basis functions are closer to the local eigenfunctions, which does not necessarily lead to smaller errors. Typically, achieving excellent results requires only two to three steps. For LKSI, the space obtained at the current iteration step is included in the space obtained at the next iteration step. Therefore, as the number of iteration steps increases, the errors will decrease, and the rate of decrease will slow down.

\section{Conclusions} \label{section:num7}

\label{sec:conclusions}
In this paper, we proposed two local subspace iteration methods for elliptic multiscale problems. Localization and the inverse operator are fundamental components of several multiscale methods. Multiple implementations of orthogonal decomposition establish the relationship between the LOD and spectral problem algorithms. Orthogonal decomposition can be regarded as an iteration step in an algorithm designed to solve spectral problems. We presented two compelling examples (the LSSI and LKSI) to illustrate our novel perspective: new multiscale methods can be designed through spectral problem algorithms. Numerical examples demonstrated that the proposed methods exhibit exceptional efficiency and applicability in  long-channel diffusion fields, which are challenging for most of   previous  multiscale methods .

This study focused  on the implementation of localization by enforcing the local homogeneous Dirichlet boundary condition. More investigation of localization would be a worthwhile research in the future. For example, in the Generalized Multiscale Finite Element Method (GMsFEM), homogeneous Neumann boundary conditions were used in the localization process. Furthermore, for asymmetric and non-positive definite operators, multiscale methods typically required  a special design \cite{Guan_Regularized_2024}. Future research will investigate asymmetric and non-positive-definite problems using the proposed multiscale methods.


\appendix

\bibliographystyle{siamplain}
\bibliography{thermoref}
\end{document}